\documentclass[11pt]{article}
\usepackage{fullpage}
\usepackage[utf8]{inputenc}
\usepackage[T1]{fontenc}  
\usepackage{amsmath}
\usepackage{amsfonts,amssymb}
\usepackage[margin=0.75in]{geometry}
\usepackage{url,hyperref}
\usepackage{amsthm,tikz}
\usepackage{xcolor}
\usepackage{float}
\usepackage{booktabs} 
\usepackage{graphicx}
\usepackage{caption}
\usepackage{subcaption}
\usepackage{ulem}
\usepackage{bm}
\allowdisplaybreaks

\newcommand{\BEAS}{\begin{eqnarray*}}
\newcommand{\EEAS}{\end{eqnarray*}}
\newcommand{\BEQ}{\begin{equation}}
\newcommand{\EEQ}{\end{equation}}
\newcommand{\BIT}{\begin{itemize}}
\newcommand{\EIT}{\end{itemize}}

\newcommand{\eg}{{\it e.g.}}
\newcommand{\ie}{{\it i.e.}}

\newcommand{\E}{\mathbb{E}}

\def\<#1,#2>{\langle #1,#2\rangle}

\newcommand{\bpi}{{\boldsymbol{\pi}}}

\newtheorem{theorem}{Theorem}
\newtheorem{remark}{Remark}
\newtheorem{assumption}{Assumption}
\newtheorem{lemma}[theorem]{Lemma}
\newtheorem{corollary}[theorem]{Corollary}

\theoremstyle{definition}
\newtheorem{definition}{Definition}[section]

{\begin{quote}}{\end{quote}}

\makeatletter
\long\def\@makecaption#1#2{
   \vskip 9pt
   \begin{small}
   \setbox\@tempboxa\hbox{{\bf #1:} #2}
   \ifdim \wd\@tempboxa > 5.5in
        \begin{center}
        \begin{minipage}[t]{5.5in}
        \addtolength{\baselineskip}{-0.95pt}
        {\bf #1:} #2 \par
        \addtolength{\baselineskip}{0.95pt}
        \end{minipage}
        \end{center}
   \else
	\hbox to\hsize{\hfil\box\@tempboxa\hfil}
   \fi
   \end{small}\par
}
\makeatother

\newcounter{oursection}

\newcounter{lecture}

\newcommand{\norm}[1]{{\left\vert\kern-0.25ex\left\vert\kern-0.25ex\left\vert #1
    \right\vert\kern-0.25ex\right\vert\kern-0.25ex\right\vert}}

\title{Homogenization and Mean-Field Approximation\\
for Multi-Player Games}

\author{Rama CONT \hfill Anran HU\\
Mathematical Institute\hfill IEOR Department\\
University of Oxford\hfill \qquad \qquad Columbia University, New York\\
{\tt Rama.Cont@maths.ox.ac.uk}\hfill {\tt ah4277@columbia.edu}}

\date{Revised version: 2026.}

\begin{document}

\maketitle

\begin{abstract}
    We investigate how the framework of mean-field games may be used to study strategic interactions in large heterogeneous populations.  Starting from a partition of the player population into 
groups, we introduce an intermediate finite-player game of mean-field type and derive explicit non-asymptotic bounds for the average exploitability of strategy profiles obtained by lifting strategies
from the associated multi-population mean-field game. The approximation error decomposes into two components: a finite-population mean-field error, controlled by empirical-measure approximation within each
group, and a heterogeneity error measuring deviations of the original players' rewards and transition dynamics from their group-level approximations. Our results apply to compact state and action spaces
and allow heterogeneous deterministic initial states within each
population. We further study the resulting trade-off between group
size and intra-group heterogeneity. 

In a parametrized heterogeneous setting, the choice of partition minimizing the resulting certified upper bound can be formulated as a mixed-integer second-order cone program. In the large-population regime, this problem is shown to be related to $K$-means clustering.
\end{abstract}
\newpage
\tableofcontents

\newpage
\section{Introduction}
  
 {\it Mean-field games} (MFGs) \cite{huang2006large,lasry2007mean} offer a tractable paradigm for the study of strategic interaction in large {\it homogeneous} populations where peer group effects are present. 
This paradigm is predicated on the assumption of a homogeneous population of agents, where the influence of any individual agent is infinitesimally small, thereby justifying the reduction of the multi-agent dynamics to the interaction between a representative agent and a ``mean field'' representing the aggregate effect of the population. 
  
The  assumptions on homogeneity of the population and symmetry of payoffs lead to greater analytical tractability and this approach has enabled a detailed theoretical study of  equilibria in such games \cite{cardaliaguet2019master,carmona2013probabilistic,cousin2011mean}.
A key theoretical result is the convergence of Nash equilibria of symmetric stochastic dynamic games involving a finite homogeneous population of players to the mean-field game limit as the player count approaches infinity \cite{possamai2023}. This convergence underpins the usefulness of MFGs in approximating the equilibria of large-scale systems with symmetric and exchangeable payoff structures.

However, the assumption of perfect homogeneity and exchangeability within the agent population is an idealization that often diverges from the inherent heterogeneity observed in real-world systems.
As most socioeconomic examples of strategic interaction relate to systems with some degree of {\it heterogeneity} across agents as well as a lack of symmetry in their interactions, it is of  interest to examine how the mean-field game framework may be useful to investigate strategic interactions and equilibria in large populations in the presence of heterogeneity. Previous work in this direction has focused on particular forms of heterogeneity which enable an approximation via a multi-population MFG with several homogeneous sub-populations \cite{bauso2016,bensoussan2018mean}, random homogeneous interaction  graphs \cite{cont2010} or graphon MFGs \cite{aurell2022,bayraktar2022stationarity,caines2019graphon,lacker2023label} on particular network structures that capture  certain departures from full symmetry in the players' interaction. 
Without such structural assumptions, the presence of heterogeneity poses significant analytical challenges and may invalidate the applicability of the mean-field paradigm to real-world systems \cite{yardim2024mean}.

A general quantitative framework for approximating heterogeneous multi-player games by MFGs, with explicit control of the error induced by different sources of heterogeneity, remains relatively underdeveloped.
Such challenges are further magnified for {\it dynamic} multi-player games. 

\paragraph{Contributions.}

We propose a systematic approach for deriving mean-field approximations for multi-player games along with {\it quantitative} error estimates. Our approach is analogous to  {homogenization} methods \cite{papanicolaou1995} used for deriving effective equations for large-scale behavior of complex heterogeneous physical systems.

Starting from a large system of players with strategic interactions, the idea is to partition the population into sub-populations and approximate the players within each group by a common group-level model.
This partition is then used to design a multi-population mean-field game \cite{bensoussan2018mean}. We provide an explicit non-asymptotic bound for the
average exploitability of strategy
profiles obtained by lifting strategies from this auxiliary
mean-field game.  Our main result, Theorem \ref{thm:homogenization}, 
shows that the average
exploitability of the lifted profile is bounded by the average
exploitability of the corresponding mean-field profile plus two
explicit error terms: a finite-population mean-field approximation
error and a heterogeneity error. The result applies to compact state
and action spaces and allows for heterogeneous deterministic initial
states within each group.

Our homogenization approach may be viewed as a natural extension of
the representative-agent approach, which corresponds to treating the
entire population as a single group. 
Here, by contrast, we allow for multiple agent types/sub-populations and, most importantly, we are able to {\it quantify} the accuracy of  mean-field approximation in terms of deviations from homogeneity and symmetry inside each sub-population.

Importantly, our approximation bounds assume no homogeneity within groups and apply
to {\it any partition} of the population. Instead, the quality of a given partition is quantified by
the heterogeneity error. The mean-field approximation error is
controlled by empirical-measure approximation within each group,
whereas the heterogeneity error measures the deviation of players' reward and transition coefficients and their
group-level approximations.

These errors are related through the process by which the population is partitioned into groups, presenting a trade-off that is pivotal to the homogenization process. 
Larger groups generally
reduce the finite-population mean-field error, whereas a finer
partition may reduce within-group heterogeneity. Using our
non-asymptotic estimates, we therefore formulate the choice of
partition as the minimization of a certified upper bound on average
exploitability. In the parametrized setting of Section~\ref{sec:partition}, this optimal partitioning problem leads to a mixed-integer
second-order cone formulation. For large populations, when the
finite-population mean-field term becomes negligible, minimization of
the heterogeneity component is closely related to the classical
$K$-means clustering problem \cite{pollard1982quantization}.

\paragraph{Related literature.} 
Mean-field approximation for homogeneous multi-player  games has been studied in various settings. In the continuous-time setting, Lasry and Lions \cite{lasry2007mean}   studied mean-field limits of {\it symmetric and homogeneous} multi-player games, where the focus is on the number of players \cite{possamai2023}. Non-asymptotic estimates were obtained by Carmona and Delarue \cite{carmona2013probabilistic} for games with affine dynamics, and more recently for more general games under weak formulations \cite{carmona2015probabilistic,lacker2016general} and with action couplings \cite{lauriere2022convergence}. These estimates were improved in \cite{delarue2019master,delarue2020master}. In the discrete-time setting, large $N$ limits were studied   in \cite{saldi2018markov,anahtarci2023,cui2021approximately}; non-asymptotic bounds have been derived  for static/stage games \cite{yardim2023stateless}, action-decoupled games with i.i.d. initial states across agents \cite{yardim2024mean}, and games with independent dynamics \cite{hu2024mf}. We refer to the recent survey in Possama\"i and Tangpi \cite{possamai2023} for further references.

Beyond the homogeneous setting, mean-field approximations have been studied for multi-population games and graph-based games, starting with the pioneering work of Huang et al. \cite{huang2006large}, which studied certain mean-field type games motivated by large communication networks. More general approximation results for multi-population games were also established in \eg, \cite{cirant2015multi,wikecek2024multiple}. 
More recently, graphon-based models have been studied in continuous-time games \cite{aurell2022,caines2019graphon,lacker2022case,lacker2023label}, discrete-time games \cite{cui2021learning}, interacting particle systems \cite{bayraktar2022stationarity}, and static games \cite{parise2023}.

Graphon models encode heterogeneity through a prescribed continuum interaction kernel \cite{parise2023}, which captures graph structures at scale $\sim N$, and finite partitions of the graphon lead naturally to multi-population games. Our construction proceeds in the opposite direction: starting from a heterogeneous finite-player game without an assumed graphon structure, we choose a partition and quantify the error incurred by replacing agents within each group by a common mean-field-type model. Thus, while graphon block approximations and our homogenization procedure share the idea of grouping agents into subpopulations, the source and role of the partition are different: in graphon models the partition approximates a given underlying interaction kernel, whereas in our framework it is an {\it optimization variable} balancing within-group heterogeneity against finite-population mean-field error.

Our results differ from existing approximation results in allowing simultaneous heterogeneity in rewards, dynamics, and deterministic initial states within a non-asymptotic multi-population homogenization framework. In addition, we study the choice of population partition as an optimization problem, balancing finite-population mean-field approximation error against within-group heterogeneity. Related work by Yardim et al.~\cite{yardim2024exploiting} studies single-population approximations for approximately symmetric multi-player games
with i.i.d.\ initial states.
\section{Mean-field games as approximations of $N$-player games}\label{background_and_setup}
Let $S\subset\mathbb R^{d_S}$ and $A\subset\mathbb R^{d_A}$ be compact sets representing respectively the {\it state space} and the set of {\it actions}. We set
\[
E:=S\times A\subset\mathbb R^d,
\qquad d=d_S+d_A,
\]
which is a compact subset of $\mathbb R^d$ equipped with the  metric
\[
d_E((s,a),(s',a')):=|s-s'|+|a-a'|.
\]
We assume, without loss of generality, $\operatorname{diam}(E)\le 1$.
Let $\mathcal P(E)$ denote the set of Borel probability measures on $E$.
For $\mu,\nu\in\mathcal P(E)$ define the Wasserstein-1 distance
\[
W_1(\mu,\nu)=\inf_{\gamma\in\Gamma(\mu,\nu)}
\int_{E\times E}|x-y|\,d\gamma(x,y)
=
\sup_{ \operatorname{Lip}(f)\le 1}
\left|
\int_E f\,d\mu-\int_E f\,d\nu
\right|
\]
where $\Gamma(\mu,\nu)$ denotes the set of couplings of $\mu$ and $\nu$, and the second equality follows from the Kantorovich duality theorem \cite{rachev1998}. 
For probability measures on $S$, we define the total variation distance by
\[
\|\mu-\nu\|_{TV}
:=
\sup_{\|f\|_\infty\le 1}
\left|
\int_S f\,d\mu-\int_S f\,d\nu
\right|.
\]

\subsection{Dynamic $N$-player games}
\label{N-player-games-setup}
We consider a finite-horizon dynamic Markov game
$$
G
=
\left(
 S, A,
\{P_t^i,R_t^i\}_{t\in\mathcal T,\;i\in[N]},
\mathbf s_0
\right)
$$
with $N$ players, where $ S$ and $ A$ are the common
compact state and action spaces introduced above,
$$
[N]:=\{1,\ldots,N\},
\qquad
\mathcal T:=\{0,\ldots,T\},
$$
and
$
\mathbf s_0=(s_0^1,\ldots,s_0^N)\in S^N
$
is the initial state profile.

For each $t\in\mathcal T$, the state and action profiles are denoted by
$$
\mathbf s_t=(s_t^1,\ldots,s_t^N)\in S^N,
\qquad
\mathbf a_t=(a_t^1,\ldots,a_t^N)\in A^N.
$$
Player $i$ receives the reward
$R_t^i(\mathbf s_t,\mathbf a_t),$
where
$
R_t^i: S^N\times A^N\longrightarrow\mathbb R
$
is Borel measurable.

For $t=0,\ldots,T-1$, the state transition of player $i$ is described
by a Borel stochastic kernel
$$
P_t^i(\,\cdot\,\vert\mathbf s,\mathbf a)
\in\mathcal P( S),
\qquad
(\mathbf s,\mathbf a)\in
 S^N\times A^N.
$$
Conditional on the current state-action profile
$(\mathbf s_t,\mathbf a_t)$, the players' next states are assumed to be
independent. Thus the conditional law of the next state profile is
$$
\mathbb P
\left(
d\mathbf s_{t+1}
\,\middle|\,
\mathbf s_t,\mathbf a_t
\right)
=
\bigotimes_{i=1}^N
P_t^i
\left(
ds_{t+1}^i
\,\middle|\,
\mathbf s_t,\mathbf a_t
\right).
$$

We restrict attention to decentralized Markov policies depending only
on a player's own current state. A policy
$
\pi=\{\pi_t\}_{t\in\mathcal T}
$
is a collection of Borel stochastic kernels
$$
\pi_t(\,\cdot\,\vert s)\in\mathcal P( A),
\qquad
s\in S.
$$
We denote by $\Pi$ the admissible class of such policies. Additional
regularity conditions on the admissible policy class will be specified below when needed.
A strategy profile is
$$
\boldsymbol\pi
=
(\pi^1,\ldots,\pi^N)\in\Pi^N.
$$
Under $\boldsymbol\pi$, conditional on the state profile $\mathbf s_t$,
the players choose their actions independently according to
$$
\mathbb P^{\boldsymbol\pi}
\left(
d\mathbf a_t
\,\middle|\,
\mathbf s_t
\right)
=
\bigotimes_{i=1}^N
\pi_t^i
\left(
da_t^i
\,\middle|\,
s_t^i
\right).
$$
Together with the transition kernels above and the deterministic
initial condition $\mathbf s_0$, these kernels uniquely determine the
law of the state-action process
$
(\mathbf s_t,\mathbf a_t)_{t\in\mathcal T}.
$

The expected cumulative reward of player $i$ under
$\boldsymbol\pi$ is
$$
V^i(\boldsymbol\pi;G)
:=
\mathbb E^{\boldsymbol\pi}
\left[
\sum_{t=0}^T
R_t^i(\mathbf s_t,\mathbf a_t)
\right].
$$
For a strategy profile $\boldsymbol{\pi}\in\Pi^N$, define the
\emph{exploitability of player $i$} by
$$
\operatorname{Expl}_i(\boldsymbol{\pi};G)
:=
\sup_{\widetilde{\pi}^i\in\Pi}
\left\{
V^i(\boldsymbol{\pi}^{-i},\widetilde{\pi}^i;G)
-
V^i(\boldsymbol{\pi};G)
\right\},
$$
where $(\boldsymbol{\pi}^{-i},\widetilde{\pi}^i)$ denotes the strategy
profile obtained by replacing $\pi^i$ by $\widetilde{\pi}^i$.

The \emph{average exploitability} of $\boldsymbol{\pi}$ is defined by
$$
\operatorname{AvgExpl}(\boldsymbol{\pi};G)
:=
\frac{1}{N}
\sum_{i=1}^N
\operatorname{Expl}_i(\boldsymbol{\pi};G).
$$

\begin{definition}[Nash equilibrium]
A strategy profile $\boldsymbol{\pi}\in\Pi^N$ is a Nash equilibrium of
$G$ if
$$
\operatorname{Expl}_i(\boldsymbol{\pi};G)=0
\qquad
\text{for every }i\in[N].
$$
Equivalently,
$
\operatorname{AvgExpl}(\boldsymbol{\pi};G)=0.
$
\end{definition}

\begin{remark}
The condition
$
\operatorname{AvgExpl}(\boldsymbol{\pi};G)\leq\epsilon
$
should be distinguished from the  $\epsilon$-Nash equilibrium
condition
$$
\max_{i\in[N]}
\operatorname{Expl}_i(\boldsymbol{\pi};G)
\leq\epsilon.
$$
The former controls the population-average gain from unilateral
deviations, whereas the latter controls the gain of every individual
player. 
\end{remark}

\subsection{Mean-field games and multi-population mean-field games}

Mean-field games \cite{huang2006large,lasry2007mean} provide tractable approximations of large symmetric
dynamic games in which individual players interact through aggregate
population distributions. To accommodate heterogeneity across groups,
we consider a multi-population mean-field game with $K$ populations associated with a partition 
$$
[N]=\bigsqcup_{k=1}^K I_k,
\qquad
N_k:=|I_k|,
\qquad
\lambda_k:=\frac{N_k}{N}>0,
\qquad
\sum_{k=1}^K\lambda_k=1.
$$
The coefficients $\lambda_k$ represent the relative sizes of the
populations.
A $K$-population mean-field game is specified by
$$
G^{\mathrm{MF}}
=
\left(
 S, A,
\{\bar P_t^k,\bar R_t^k\}_{t\in\mathcal T,\;k\in[K]},
\{\mu_0^k\}_{k\in[K]}
\right),
$$
where $ S$ and $ A$ are the compact state and action
spaces introduced above, and
$
\mu_0^k\in\mathcal P( S)
$
is the initial state distribution of population $k$.

For each $t\in\mathcal T$ and $k\in[K]$, the reward function is a
Borel measurable map
$$
\bar R_t^k:
 S\times A\times\mathcal P(E)^K
\longrightarrow
\mathbb R.
$$
For $t=0,\ldots,T-1$, the transition kernel is a Borel stochastic
kernel
$$
\bar P_t^k
\left(
\,\cdot\,
\middle|
s,a,L^1,\ldots,L^K
\right)
\in\mathcal P( S),
\qquad
{\rm for}\quad
(s,a,L^1,\ldots,L^K)
\in
 S\times A\times\mathcal P(E)^K.
$$

A mean-field strategy profile is
$$
\bar{\boldsymbol\pi}
=
(\bar\pi^1,\ldots,\bar\pi^K)
\in\Pi^K,
$$
where $\bar\pi^k$ is the policy used by population $k$.

Given $\bar{\boldsymbol\pi}$, let
$\mu_t^{k,\bar{\boldsymbol\pi}}
\in\mathcal P( S)
$
denote the state distribution of population $k$ at time $t$, and let
$
L_t^{k,\bar{\boldsymbol\pi}}
\in\mathcal P(E)
$
denote its joint state-action distribution. These measures are
defined recursively by
$$
L_t^{k,\bar{\boldsymbol\pi}}(ds,da)
=
\mu_t^{k,\bar{\boldsymbol\pi}}(ds)
\bar\pi_t^k(da\mid s),
$$
with
$
\mu_0^{k,\bar{\boldsymbol\pi}}=\mu_0^k,
$
and, for $t=0,\ldots,T-1$,
$$
\mu_{t+1}^{k,\bar{\boldsymbol\pi}}(B)
=
\int_{ S\times A}
\bar P_t^k
\left(
B
\middle|
s,a,
L_t^{1,\bar{\boldsymbol\pi}},
\ldots,
L_t^{K,\bar{\boldsymbol\pi}}
\right)
L_t^{k,\bar{\boldsymbol\pi}}(ds,da)
$$
for every Borel set $B\subseteq S$.
The corresponding joint state-action distribution is
$$
L_{t+1}^{k,\bar{\boldsymbol\pi}}(ds',da')
=
\mu_{t+1}^{k,\bar{\boldsymbol\pi}}(ds')
\bar\pi_{t+1}^k(da'\mid s').
$$
Thus the strategy profile $\bar{\boldsymbol\pi}$ determines a
flow
$
\left(
L_t^{1,\bar{\boldsymbol\pi}},
\ldots,
L_t^{K,\bar{\boldsymbol\pi}}
\right)_{t\in\mathcal T}
$
of population state-action distributions.

To define unilateral deviations, fix a population $k$ and a policy
$\widetilde\pi^k\in\Pi$. A representative player from population $k$
uses the policy $\widetilde\pi^k$, while the population distributions
remain those generated by the reference profile
$\bar{\boldsymbol\pi}$. More precisely, the representative player's
initial state has law $\mu_0^k$, its action at time $t$ is sampled
according to
$$
\widetilde\pi_t^k(da\mid s),
$$
and its state evolves according to
$$
\bar P_t^k
\left(
ds'
\middle|
s,a,
L_t^{1,\bar{\boldsymbol\pi}},
\ldots,
L_t^{K,\bar{\boldsymbol\pi}}
\right).
$$

The corresponding expected cumulative reward is denoted by
$$
\bar V^k
\left(
\widetilde\pi^k,\bar{\boldsymbol\pi};
G^{\mathrm{MF}}
\right)
:=
\mathbb E
\left[
\sum_{t=0}^T
\bar R_t^k
\left(
s_t^k,a_t^k,
L_t^{1,\bar{\boldsymbol\pi}},
\ldots,
L_t^{K,\bar{\boldsymbol\pi}}
\right)
\right].
$$
To account for possibly heterogeneous initial states within each
population, it is useful to make the dependence of the representative
player's value on its initial state explicit. For $s\in S$,
we denote by
$$
\bar V^{k,s}
\left(
\widetilde\pi^k,\bar{\boldsymbol\pi};
G^{\mathrm{MF}}
\right)
$$
the expected cumulative reward of a representative player from
population $k$ which starts from the deterministic initial state
$s$, uses the policy $\widetilde\pi^k$, and evolves in the
mean-field environment
$$
\left(
L_t^{1,\bar{\boldsymbol\pi}},
\ldots,
L_t^{K,\bar{\boldsymbol\pi}}
\right)_{t\in\mathcal T}
$$
generated by the reference profile $\bar{\boldsymbol\pi}$.
More precisely,
$$
\bar V^{k,s}
\left(
\widetilde\pi^k,\bar{\boldsymbol\pi};
G^{\mathrm{MF}}
\right)
:=
\mathbb E
\left[
\sum_{t=0}^T
\bar R_t^k
\left(
s^k_t,a^k_t,
L_t^{1,\bar{\boldsymbol\pi}},
\ldots,
L_t^{K,\bar{\boldsymbol\pi}}
\right)
\right],
$$
where $s^k_0=s$, the actions are sampled according to
$\widetilde\pi_t^k(\cdot\mid s_t)$, and, for
$t=0,\ldots,T-1$,
$$
s^k_{t+1}
\sim
\bar P_t^k
\left(
\cdot
\mid
s^k_t,a^k_t,
L_t^{1,\bar{\boldsymbol\pi}},
\ldots,
L_t^{K,\bar{\boldsymbol\pi}}
\right).
$$

The \emph{state-conditioned exploitability} of population $k$ at
initial state $s$ is defined by
$$
\operatorname{Expl}^{k,s}
\left(
\bar{\boldsymbol\pi};
G^{\mathrm{MF}}
\right)
:=
\sup_{\widetilde\pi^k\in\Pi}
\left\{
\bar V^{k,s}
\left(
\widetilde\pi^k,\bar{\boldsymbol\pi};
G^{\mathrm{MF}}
\right)
-
\bar V^{k,s}
\left(
\bar\pi^k,\bar{\boldsymbol\pi};
G^{\mathrm{MF}}
\right)
\right\}.
$$

We define the exploitability of population $k$ by averaging the
state-conditioned exploitability with respect to its initial
distribution:
$$
\operatorname{Expl}_k
\left(
\bar{\boldsymbol\pi};
G^{\mathrm{MF}}
\right)
:=
\int_{ S}
\operatorname{Expl}^{k,s}
\left(
\bar{\boldsymbol\pi};
G^{\mathrm{MF}}
\right)
\mu_0^k(ds).
$$
The average exploitability of
$\bar{\boldsymbol\pi}$ is then
$$
\operatorname{AvgExpl}
\left(
\bar{\boldsymbol\pi};
G^{\mathrm{MF}}
\right)
:=
\sum_{k=1}^K
\lambda_k
\operatorname{Expl}_k
\left(
\bar{\boldsymbol\pi};
G^{\mathrm{MF}}
\right).
$$

In particular, when the initial distribution of population $k$ is
the empirical distribution
$$
\mu_0^k
=
\mu_0^{k,N}
=
\frac{1}{N_k}
\sum_{i\in I_k}
\delta_{s_0^i},
$$
we have
$$
\operatorname{Expl}_k
\left(
\bar{\boldsymbol\pi};
G^{\mathrm{MF}}
\right)
=
\frac{1}{N_k}
\sum_{i\in I_k}
\operatorname{Expl}^{k,s_0^i}
\left(
\bar{\boldsymbol\pi};
G^{\mathrm{MF}}
\right),
$$
and therefore
$$
\operatorname{AvgExpl}
\left(
\bar{\boldsymbol\pi};
G^{\mathrm{MF}}
\right)
=
\frac{1}{N}
\sum_{k=1}^K
\sum_{i\in I_k}
\operatorname{Expl}^{k,s_0^i}
\left(
\bar{\boldsymbol\pi};
G^{\mathrm{MF}}
\right).
$$

\begin{definition}[Nash equilibrium for a multi-population mean-field game]
A mean-field strategy profile
$\bar{\boldsymbol\pi}\in\Pi^K$ is a Nash equilibrium of
$G^{\mathrm{MF}}$ if
$$
\operatorname{Expl}^{k,s}
\left(
\bar{\boldsymbol\pi};G^{\mathrm{MF}}
\right)
=0
$$
for $\mu_0^k$-almost every $s\in S$ and every
$k\in[K]$. Equivalently,
$
\operatorname{AvgExpl}
\left(
\bar{\boldsymbol\pi};G^{\mathrm{MF}}
\right)
=0.
$
\end{definition}

For $\epsilon\geq0$, we say that the mean-field strategy profile
$\bar{\boldsymbol\pi}$ has \emph{average exploitability
at most $\epsilon$} if
$$
\operatorname{AvgExpl}
\left(
\bar{\boldsymbol\pi};G^{\mathrm{MF}}
\right)
\leq\epsilon.
$$

\begin{remark}
Players first evaluate the gain from an optimal unilateral
deviation conditional on the representative player's initial state
and then we average these gains over the population
distribution. Thus
$$
\operatorname{Expl}_k
\left(
\bar{\boldsymbol\pi};G^{\mathrm{MF}}
\right)
=
\int_{ S}
\sup_{\widetilde\pi^k\in\Pi}
\left\{
\bar V^{k,s}
\left(
\widetilde\pi^k,\bar{\boldsymbol\pi};
G^{\mathrm{MF}}
\right)
-
\bar V^{k,s}
\left(
\bar\pi^k,\bar{\boldsymbol\pi};
G^{\mathrm{MF}}
\right)
\right\}
\mu_0^k(ds),
$$
rather than placing the supremum outside the integral. This
criterion is the natural mean-field counterpart of average
exploitability in the finite-player game when players have
heterogeneous deterministic initial states. If there exists a
single admissible Markov policy which is simultaneously optimal for
$\mu_0^k$-almost every initial state, the two formulations coincide.
\end{remark}

Let us motivate the problem and our methodology using an  example adapted from \cite{cousin2011mean,guo2023general}.

\subsection{Example: Equilibrium pricing}\label{subsec:example}

In a market scenario involving a large number of firms, each producing the same product from some raw material, strategic decision-making plays a crucial role. We examine an $N$-player game over some finite horizon, where each of the $N$ firms independently decides on production quantity $q_t^i$ and inventory replenishment $h_t^i$ of the raw material at each time $t$. The current inventory level of raw material of firm $i$ at time $t$, denoted as $s_t^i$, evolves according to the equation:
\[
s_{t+1}^i = s_t^i - \min\{q_t^i, s_t^i\} + h_t^i.
\]
In the context of the general notation of $N$-player games in Section \ref{N-player-games-setup}, the states are $s_t^i$, the actions are $a_t^i=(q_t^i,h_t^i)$, and the above equation defines the (deterministic) transition model $P_t^i$ ($i\in[N]$). To fit into the framework in this section, we set $S=\{0,1,\ldots,S_{\max}\}$ and $A=\{(q,h):q\in\{0,1,\ldots,Q\},\;
h\in\{0,1,\ldots,H\}\}$ for some positive integers $S_{\max}$, $Q$ and $H$.

We assume the firm will pay a cost for an emergency order of the raw material if the firm overproduces and exceeds its current inventory capacity (\ie, $q_t^i>s_t^i$). The (deterministic) reward for firm $i$ at time $t$ is given by:
\[
r_t^i({\bf s}_t,{\bf a}_t)=R_t^i({\bf s}_t,{\bf a}_t)= (\hat{p}_t - c_0^i)q_t^i - c_1^i(q_t^i)^2 - c_2^ih_t^i - (c_2^i + c_3^i)\max\{q_t^i - s_t^i, 0\} - c_4^is_t^i,
\]
where $\hat{p}_t$ represents the market price of the product at time $t$, $c_0^i$ covers per-unit manufacturing and labor costs, $c_1^i$ denotes a quadratic transient price impact cost associated with the production level, $c_2^i$ is the cost for regular raw material orders, $c_3^i$ represents the additional cost for emergency raw material orders, and $c_4^i$ is the inventory holding cost. We remark that these costs are usually different across the firms. For example, different firms may employ  different technologies for production and inventory management. They may also have different supply chain strategies, different relationships with suppliers, and different geographical locations. These factors can lead to different cost coefficients.


The price  $\hat{p}_t$ is  endogenously determined by a supply-demand equilibrium \cite{bernstein2006regional}. We allow for a fixed baseline supply $NQ_0$, with $Q_0>0$, which is exogenous and is not an optimization variable in the present model.
The total supply is therefore
$NQ_0+\sum_{i=1}^N q_t^i$, while the total demand is 
$d_N \hat{p}_t^{-\sigma}$, where $d_N = \delta\ N$ reflects the linear scaling of demand with the number of producers, implying that the number of customers grows proportionally to the number of firms. Here $\delta $ denotes some benchmark unit demand level and $\sigma>0$ is the elasticity of demand that can be
interpreted as the elasticity of substitution between the given product and any other good. 
This results in the equilibrium condition: 
\[
Q_0+\frac{1}{N}\sum_{i=1}^N q_t^i = \frac{\delta}{ \hat{p}_t^{\sigma}}.
\]
The baseline supply $Q_0$ may be interpreted as an exogenous policy intervention, while the positivity of $Q_0$ also prevents the inverse demand function from becoming singular when aggregate production is zero.

In this example, if all the coefficients $(c_m^i),m\in\{0,1,2,3,4\}$ are the same across all firms, as $N$ approaches infinity, the above model admits a natural mean-field counterpart. Transitioning from the $N$-player setting to the mean-field framework, we consider a continuum of homogeneous firms, simplifying the analysis by focusing on a representative firm's interactions with the population distribution. In this mean-field setting, each firm's decisions on production $q_t$ and inventory replenishment $h_t$ are based on the representative firm's strategy, with the inventory level $s_t$ evolving similarly to the individual firm case in the $N$-player game. The reward function in the mean-field setting mirrors that of the $N$-player game, with the market price $p_t$ now determined by the mean-field supply-demand equilibrium: 
\[
Q_0+\E[q_t]=\frac{\delta}{p_t^{\sigma}},
\]
where $\E[q_t]$ is the averaged supply of all firms, and the expectation is over the $q$-marginal of the (mean-field) joint state-action distribution. We refer  to \cite{cousin2011mean,guo2023general} for more details. 

When the coefficients $(c_m^i),m\in\{0,1,2,3,4\}$ are not the same across  firms,  a representative agent model is only an approximation whose accuracy needs to be examined. We shall return to this example below.

The present example focuses on heterogeneity in rewards while keeping the inventory dynamics homogeneous across firms. Heterogeneous transition dynamics can be treated within the general framework developed below.

\section{Homogenization of dynamic multi-player games} \label{sec:homogenization}
We now describe a general procedure for approximating a multi-player dynamic game by a (multi-population) mean-field game.

Starting from a general multi-player game with heterogeneous interactions, agents with similar characteristics can be grouped together
into {\it near-homogeneous} groups, which we can then approximate as being homogeneous.  
We can then replace the heterogeneous interactions between agents in such clusters with interactions between representative agents of each cluster. 
This is the {\it symmetrization} step.

This leads to an approximation of the original game by a {\it multi-population game of mean-field type}, in which agents within each group share common reward and transition models and their interactions can be described in terms of a representative agent of each group and the corresponding joint state-action distributions, facilitating a transition to a multi-population mean-field game (MFG) framework. 

The complete homogenization process consists of two steps. The first step is to partition the $N$ players into some disjoint subsets and homogenize each group, i.e. replace it by identical agents with common reward function and transition kernel. The second step is to approximate the behavior of each group by its large population limit.  
To formalize this procedure, we introduce an intermediate model that serves as the bridge between the original $N$ player game $G$ and the limiting MFG $G^{\mathrm{MF}}$.

\paragraph{$N$-player games of mean-field type.} 
We next introduce a class of $N$-player games with a mean-field
interaction structure. This class will serve as an intermediate model
between the original heterogeneous $N$-player game $G$ and the
multi-population mean-field game $G^{\mathrm{MF}}$.
Consider a partition of the population into $K$  groups
$$
[N]=\bigsqcup_{k=1}^K I_k,
\qquad
N_k:=|I_k|.
$$
For a state-action profile
$(\mathbf s,\mathbf a)\in S^N\times A^N$, define the
empirical state-action distribution of group $k$ by
$$
L^{k,N}(\mathbf s,\mathbf a)
:=
\frac{1}{N_k}
\sum_{i\in I_k}
\delta_{(s^i,a^i)}
\in\mathcal P(E).
$$
Along the state-action process, we write
$$
L_t^{k,N}
:=
\frac{1}{N_k}
\sum_{i\in I_k}
\delta_{(s_t^i,a_t^i)}.
$$

We say that an $N$-player game is of \emph{mean-field type} with
respect to the partition $(I_k)_{k=1}^K$ if there exist reward
functions
$$
\bar R_t^k:
 S\times A\times\mathcal P(E)^K
\longrightarrow\mathbb R
$$
and stochastic kernels
$$
\bar P_t^k
\left(
\,\cdot\,
\middle|
s,a,L^1,\ldots,L^K
\right)
\in\mathcal P( S)
$$
such that, for every player $i\in I_k$,
$$
\widehat R_t^i(\mathbf s,\mathbf a)
=
\bar R_t^k
\left(
s^i,a^i,
L^{1,N}(\mathbf s,\mathbf a),
\ldots,
L^{K,N}(\mathbf s,\mathbf a)
\right),
$$
and, for $t=0,\ldots,T-1$,
$$
\widehat P_t^i
\left(
ds'
\middle|
\mathbf s,\mathbf a
\right)
=
\bar P_t^k
\left(
ds'
\middle|
s^i,a^i,
L^{1,N}(\mathbf s,\mathbf a),
\ldots,
L^{K,N}(\mathbf s,\mathbf a)
\right).
$$

We denote the resulting $N$-player game by
$G_N^{\mathrm{MF}}
:=
\left(
S,A,
\{\widehat P_t^i,\widehat R_t^i\}_{t\in\mathcal T,i\in[N]},
\mathbf s_0
\right)$. 
Thus players belonging to the same group have the same reward and
transition coefficients, and their interaction with the other players
depends only on the empirical state-action distributions
$$
\left(
L_t^{1,N},\ldots,L_t^{K,N}
\right).
$$
The players within a group may nevertheless have different
deterministic initial states.

Associated with $G_N^{\mathrm{MF}}$ is the $K$-population mean-field
game $G^{\mathrm{MF}}$ introduced in Section~2.2, with the same reward
functions $\bar R_t^k$ and transition kernels $\bar P_t^k$. Its initial
state distribution in population $k$ is chosen to be the empirical
distribution of the initial states of the corresponding finite
population:
$$
\mu_0^{k,N}
:=
\frac{1}{N_k}
\sum_{i\in I_k}
\delta_{s_0^i}
\in\mathcal P( S).
$$
Thus the mean-field model retains the empirical distribution of the
initial states within each group while discarding their labels. 

\begin{figure}[H]
 \centering
 \includegraphics[width=10cm, trim={1cm 1cm 6cm 1cm},clip]{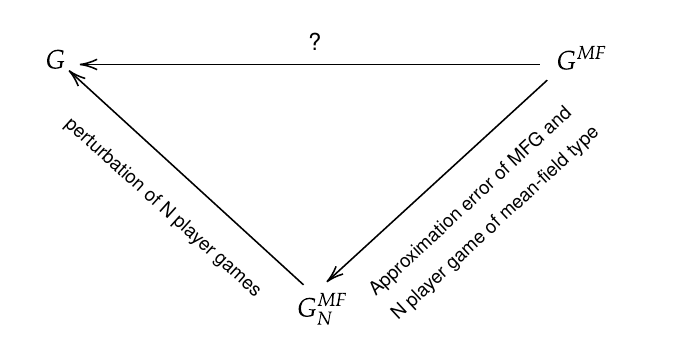}
 \caption{Decomposition of homogenization error.}\label{fig:triangle}
 \end{figure}

Given a mean-field strategy profile
$\bar{\boldsymbol\pi}=(\bar\pi^1,\ldots,\bar\pi^K)$, we lift it to an
$N$-player strategy profile $\boldsymbol\pi$ by setting
$$
\pi^i=\bar\pi^k,
\qquad i\in I_k.
$$
Our objective is to quantify the average exploitability of this lifted
profile in the original game $G$. The resulting error naturally
decomposes into two components: the error incurred by replacing
$G_N^{\mathrm{MF}}$ with its multi-population mean-field limit
$G^{\mathrm{MF}}$, and the heterogeneity error incurred by replacing
the original game $G$ with the homogenized game
$G_N^{\mathrm{MF}}$, as illustrated in Figure~\ref{fig:triangle}. We quantify these two error components separately below.



We now revisit the equilibrium pricing example of Section \ref{subsec:example}, and demonstrate how the mean-field approximation is carried out in this problem. 

\paragraph{Example revisited.} For the equilibrium pricing problem discussed in Section \ref{subsec:example},   heterogeneity arises from the cost/reward functions. Therefore the homogenization procedure can be viewed as   partitioning the firms according to their cost functions (or the coefficients defined in these functions). For example,  the set of coefficients may be partitioned into $K$ disjoint subsets $C^k\subset \mathbb{R}^5$ ($k\in[K]$); firm $i$  belongs to cluster $k$ if and only if $(c_0^i,\dots,c_4^i)\in C^k$. We then average the cost coefficients in each cluster and use a cost function with coefficients $(\bar c_0^k,\dots,\bar c_4^k)$ to represent all firms in group $k$. Using this partition of the firms $ I_k$ $(k=1,\dots,K)$ with $N_k=| I_k|$, we can construct an $N$-player game of mean-field type $G_N^{\mathrm{MF}}$ as follows. 
\begin{itemize}
    \item For firm $i\in I_k$, the reward function is given by
\begin{align*}
\hat{R}_t^i({\bf s}_t,{\bf a}_t)&= \bar{R}_t^k(s_t^i,a_t^i,L_t^{1,N},\dots,L_t^{K,N}) \\
&= (\hat p_t - \bar c_0^k)q_t^i - \bar c_1^k(q_t^i)^2 - \bar c_2^kh_t^i - (\bar c_2^k + \bar c_3^k)\max\{q_t^i - s_t^i, 0\} - \bar c_4^ks_t^i,
\end{align*}
when the firm with current inventory level  $s_t^i$  decides to produce $q_t^i$ amount of product and replenish $h_t^i$ amount of raw material. Here $\hat p_t$ solves $Q_0+\frac{1}{N}\sum_{i=1}^N q_t^i =  \hat p_t^{-\sigma}\delta $. Recall that ${\bf s}_t=(s_t^1,\dots,s_t^N)$, ${\bf a}_t=(a_t^1,\dots,a_t^N)$, with $a_t^i=(q_t^i,h_t^i)$ for $i\in[N]$.
\item The state of firm $i\in I_k$ is initialized at $s_0^i$ and evolves according to
\[
s_{t+1}^i = s_t^i - \min\{q_t^i, s_t^i\} + h_t^i.
\]
\end{itemize}
The above game is of mean-field type because one can rewrite $\frac{1}{N}\sum_{i=1}^N q_t^i$ as 
\[
\frac{1}{N}\sum_{i=1}^N q_t^i= \frac{1}{N}\sum_{k=1}^K\sum_{i\in I_k} \sum_{j=0}^Qj{\bf 1}_{q_t^i=j}=
\sum_{k=1}^K \frac{N_k}{N}\cdot \frac{1}{N_k}\sum_{i\in I_k} \sum_{j=0}^Q j{\bf 1}_{q_t^i=j} =\sum_{k=1}^K \frac{N_k}{N} \sum_{j=0}^Q j\sum_{s=0}^{S_{\max}}\sum_{h=0}^HL_t^{k,N}(s,j,h),
\]
which is a function of the joint empirical distribution of the states and actions across all firms in each group $k$: $L_t^{k,N}(s,a)=\frac{1}{N_k}\sum_{i\in I_k}{\bf 1}\{s_t^i=s,q_t^i=j,h_t^i=h\}$, where $a=(j,h)$. With slight abuse of notation, we also denote $L_t^{k,N}(s,a)$ as $L_t^{k,N}(s,j,h)$. 

Based on this $G_N^{\mathrm{MF}}$, we define the limiting $K$-population MFG as follows. 
\begin{itemize}
    \item Consider $K$ groups (continuum) of homogeneous firms, with a representative firm in each group $k$ whose reward function is
\[
\bar{R}_t^k(s^k_t,a^k_t,L_t^1,\dots,L_t^K) = (p_t - \bar c_0^k)q_t^k - \bar c_1^k(q_t^k)^2 - \bar c_2^kh_t^k - (\bar c_2^k + \bar c_3^k)\max\{q_t^k - s_t^k, 0\} - \bar c_4^ks_t^k,
\]
when the representative firm with current inventory level $s_t^k$ decides to produce $q_t^k$ amount of product and replenish $h_t^k$ amount of raw material. Here the market price $p_t$ is defined as the solution of 
\[
Q_0+\sum_{k\in[K]} \frac{N_k}{N}\sum_{q=0}^Qq\sum_{s=0}^{S_{\max}}\sum_{h=0}^HL_t^{k}(s,q,h) = \frac{\delta }{p_t^{\sigma}},
\]
with $L_t^k$ being the 
(mean-field) joint distribution of the state $s_t$ and action $a_t=(q_t,h_t)$ of all the firms in group $k$. Again, with slight abuse of notation, here $L_t^k(s,q,h)=L_t^k(s,a)$ with $a=(q,h)$. 
\item The initial state of the representative firm in group $k$ satisfies
$s_0^k\sim\mu_0^k=\mu_0^{k,N}
=
\frac{1}{N_k}\sum_{i\in I_k}\delta_{s_0^i}$, 
and its state evolves according to
\[
s_{t+1}^k = s_t^k - \min\{q_t^k, s_t^k\} + h_t^k.
\]
\item If the firms in group $k\in[K]$ adopt policy $\bar\pi^k:\mathcal T\times S\rightarrow\mathcal P(A)$, then one can explicitly characterize $L_t^k(s,q,h)=\mu_t^k(s)\bar\pi_t^k(q,h\mid s)$, where $\mu_t^k$ is the distribution of the inventory level of all firms in group $k$ which evolves according to  \[
\mu_{t+1}^k(s')=\sum_{\{s,q,h:s-\min\{q,s\}+h=s'\}}\mu_t^k(s)\bar\pi_t^k(q,h\mid s).
\]
\end{itemize}

\section{Quantitative estimates of the homogenization error}\label{sec:homogenization_error}

We now quantify the approximation error introduced by the two-step
homogenization procedure of Section~3.

Let $\bar{\boldsymbol\pi}$ be a strategy profile for the
multi-population mean-field game $G^{\mathrm{MF}}$, and let
$\boldsymbol\pi$ denote its lift to the $N$-player game, defined by
$$
\pi^i=\bar\pi^k,
\qquad i\in I_k.
$$
Our objective is to bound the average exploitability of
$\boldsymbol\pi$ for the  heterogeneous $N$-player game $G$.

\subsection{Assumptions and useful lemmas} 
We introduce the following assumptions:

\begin{assumption}[Bounded rewards]
There exists $\bar R_{\max}<\infty$ such that
\[
|\bar R_t^k(s,a,L^1,\ldots,L^K)|
\le
\bar R_{\max}
\]
for all $t,k,s,a,L^1,\ldots,L^K$.\label{ass.boundedreward}
\end{assumption}

\begin{assumption}[Weighted Lipschitz continuity in the mean field]\label{ass.weightedLipschitz}
There exist nonnegative constants
\[
w_R^1,\ldots,w_R^K,
\qquad
w_P^1,\ldots,w_P^K
\]
such that, for all $t,k,s,a$ and all
$L^1,\ldots,L^K,\widetilde L^1,\ldots,\widetilde L^K\in\mathcal P(E)$,
\[
\left|
\bar R_t^k(s,a,L^1,\ldots,L^K)
-
\bar R_t^k(s,a,\widetilde L^1,\ldots,\widetilde L^K)
\right|
\le
\sum_{j=1}^K
w_R^j
W_1(L^j,\widetilde L^j),
\]
and
\[
\left\|
\bar P_t^k(\cdot\mid s,a,L^1,\ldots,L^K)
-
\bar P_t^k(\cdot\mid s,a,\widetilde L^1,\ldots,\widetilde L^K)
\right\|_{TV}
\le
\sum_{j=1}^K
w_P^j
W_1(L^j,\widetilde L^j).
\]
\end{assumption}
We define
\begin{equation}
    W_{\max}
:=
\max\left\{
\sum_{j=1}^K w_R^j,
\sum_{j=1}^K w_P^j
\right\}.
\end{equation}

\begin{assumption}[Wasserstein continuity in own state and action]
There exists $L_P<\infty$ such that, for every
$t=0,\ldots,T-1$, $k\in[K]$, $s,s'\in S$,
$a,a'\in A$, and
$(L^1,\ldots,L^K)\in\mathcal P(E)^K$,
$$
W_1\left(
\bar P_t^k
\left(
\cdot\mid s,a,L^1,\ldots,L^K
\right),
\bar P_t^k
\left(
\cdot\mid s',a',L^1,\ldots,L^K
\right)
\right)
\leq
L_P\left(
|s-s'|+|a-a'|
\right).
$$
\label{ass.transition_lipschitz}
\end{assumption}

\begin{assumption}[Lipschitz continuity in own state and action]\label{ass.reward_lipschitz}

There exists $L_R<\infty$ such that
\[
\begin{aligned}
|\bar R_t^k(s,a,L^1,\ldots,L^K)
-
\bar R_t^k(s',a',L^1,\ldots,L^K)|
&\le
L_R(|s-s'|+|a-a'|).
\end{aligned}
\]
\end{assumption} 

The following assumption controls the
stability of policies under perturbations of states:
\begin{assumption}[Stability of policies under perturbations of states]
There exists $L_\pi<\infty$ such that, for every $\pi\in\Pi$,
$t\in\mathcal T$, and $s,s'\in S$,
$$
W_1
\left(
\pi_t(\cdot\mid s),
\pi_t(\cdot\mid s')
\right)
\leq
L_\pi |s-s'|.
$$ \label{ass.regularkernel}
\end{assumption}

\begin{remark}
For finite state spaces, Assumption~\ref{ass.regularkernel} is automatically satisfied. More generally, the Lipschitz regularity imposed on randomized policies is mainly a technical condition used to obtain the stability estimates in the analysis, and is also natural in learning and approximation settings. This restriction entails only a controlled loss of generality: under Assumptions~\ref{ass.boundedreward}, \ref{ass.transition_lipschitz}, and~\ref{ass.reward_lipschitz}, the following lemma shows that restricting deviations to policies with Lipschitz constant at most $L_\pi$ decreases the exploitability by at most $C/L_{\pi}$. Thus, the unrestricted game is recovered asymptotically as $L_\pi\to\infty$.
    \end{remark}

Let $\Pi_\infty$ denote the class of all Borel measurable policies. For any policy class $\mathcal C\subseteq\Pi_\infty$, define 
\[
\operatorname{Expl}^{k,s}_{\mathcal C}
\bigl(\bar{\bm\pi};G^{\rm MF}\bigr)
:=
\sup_{\widetilde\pi^k\in\mathcal C}
\left\{
\bar V^{k,s}
\bigl(\widetilde\pi^k,\bar{\bm\pi};G^{\rm MF}\bigr)
-
\bar V^{k,s}
\bigl(\bar\pi^k,\bar{\bm\pi};G^{\rm MF}\bigr)
\right\},
\]
and
\[
\operatorname{AvgExpl}_{\mathcal C}
\bigl(\bar{\bm\pi};G^{\rm MF}\bigr)
:=
\sum_{k=1}^K\lambda_k
\int_S
\operatorname{Expl}^{k,s}_{\mathcal C}
\bigl(\bar{\bm\pi};G^{\rm MF}\bigr)
\,\mu_0^k(ds).
\]

\begin{lemma}\label{lem:lipschitz_policy_approximation}
    Let $E= S\times  A\subset\mathbb{R}^d$ be a  compact set. For $L_\pi>0$, define 
\[
\Pi_{L_\pi}:=\left\{\pi\in\Pi_\infty:  W_1
\bigl(\pi_t(\cdot\mid s), \pi_t(\cdot\mid s')\bigr)
\le L_{\pi} |s-s'|,
\quad
t\in\mathcal T,\ s,s'\in  S
\right\}.
\]
Under Assumptions~\ref{ass.boundedreward}, \ref{ass.transition_lipschitz} and~\ref{ass.reward_lipschitz} there exists a constant $C=C( S, T, L_P, L_R)$ such that for any $k\in[K]$, every mean-field strategy profile $\bar{\bm\pi}$, and any $s\in  S$, 
\[
0\le \sup_{\widetilde\pi^k\in\Pi_\infty}
\bar V^{k,s} \bigl(
\widetilde\pi^k,\bar{\bm\pi};G^{\mathrm{MF}}
\bigr) -
\sup_{\widetilde\pi^k\in\Pi_{L_\pi}}
\bar V^{k,s}
\bigl(\widetilde\pi^k,\bar{\bm\pi};G^{\mathrm{MF}}\bigr) \le \frac{C}{L_{\pi}}.
\]
Consequently,
\[
\operatorname{AvgExpl}_{\Pi_\infty}
\bigl(\bar{\bm\pi};G^{\rm MF}\bigr)
\le
\operatorname{AvgExpl}_{\Pi_{L_\pi}}
\bigl(\bar{\bm\pi};G^{\rm MF}\bigr)
+\frac{C}{L_\pi}.
\]
\end{lemma}

The following lemma, which is a consequence of Bayraktar and Wu \cite[Lemma A.1]{bayraktar2022stationarity}, will be useful for our estimates.
\begin{lemma}[Empirical-measure estimate] 
\label{lem:empirical-BL}
Let $E\subset\mathbb{R}^d$ be compact and let $Z_1,\ldots,Z_m$ be independent, not necessarily identically distributed, $E$-valued random variables with laws $\nu_1,\ldots,\nu_m$. Define 
\[
\widehat{\nu}_m := \frac{1}{m}\sum_{\ell=1}^m\delta_{Z_\ell}, \qquad \bar{\nu}_m := \frac{1}{m}\sum_{\ell=1}^m\nu_\ell. 
\]
Then 
\[
\mathbb{E}\bigl[ W_1(\widehat{\nu}_m,\bar{\nu}_m) \bigr] \le \alpha_m(E), 
\]
where, for some constant $C_E\ge 1$ depending only on $E$, 
\begin{eqnarray}
\alpha_m(E) := C_E 
\begin{cases} 
m^{-1/2}, & d=1,\\[1mm]
m^{-1/2}\log(1+m), & d=2,\\[1mm] \label{eq.alpham}
m^{-1/d}, & d\ge 3. 
\end{cases} 
\end{eqnarray}
If $E$ is finite, one may instead take 
\[ 
\alpha_m(E)=\sqrt{\frac{|E|}{m}}. 
\]
\end{lemma}

For the explicit scaling rate examples below we will assume $d\geq 3$. The cases $d=1$ and $d=2$ follow from
Lemma~\ref{lem:empirical-BL} with the corresponding sharper rates,
with a logarithmic correction when $d=2$.

\subsection{Homogenization bound}

We now state our main result, which gives a non-asymptotic bound on the average exploitability of a strategy profile obtained by lifting a multi-population mean-field strategy to the original heterogeneous $N$-player game. The bound separates the mean-field approximation error from the error due to intra-group heterogeneity.
 \begin{theorem}[Homogenization bound]
\label{thm:homogenization}
Let $G=
\left(
S,A,
\{P_t^i,R_t^i\}_{t\in\mathcal T,i\in[N]},
\mathbf s_0
\right)$ be the original $N$-player game.
Consider a partition
\[
[N]=\bigsqcup_{k=1}^K I_k,
\qquad N_k:=|I_k|, 
\]
and let
$G_N^{\mathrm{MF}}
=
\left(
S,A,
\{\widehat P_t^i,\widehat R_t^i\}_{t\in\mathcal T,i\in[N]},
\mathbf s_0
\right)$ be an $N$-player game of mean-field type with respect to this partition.
Let $G^{\mathrm{MF}}
=
\left(
S,A,
\{\bar P_t^k,\bar R_t^k\}_{t\in\mathcal T,k\in[K]},
\{\mu_0^k\}_{k\in[K]}
\right)$ be the corresponding $K$-population
mean-field game, where $\mu_0^k
=
\frac1{N_k}\sum_{i\in I_k}\delta_{s_0^i}$. Let 
$\bar{\boldsymbol\pi}
=(\bar\pi^1,\ldots,\bar\pi^K)\in\Pi^K$
be a mean-field strategy profile, and define its lift
$\boldsymbol\pi\in\Pi^N$ by
$ \pi^i=\bar\pi^k, i\in I_k.$
Under Assumptions~\ref{ass.boundedreward},  \ref{ass.weightedLipschitz},  \ref{ass.transition_lipschitz}, \ref{ass.reward_lipschitz}, and~\ref{ass.regularkernel}, we have
$$
\operatorname{AvgExpl}(\boldsymbol\pi;G)
\leq
\operatorname{AvgExpl}
(\bar{\boldsymbol\pi};G^{\mathrm{MF}})
+
\epsilon_{\mathrm{MF}}
+
\epsilon_{\mathrm{heter}},
$$
where 
\begin{equation}
\epsilon_{\mathrm{MF}}
=
C_{\mathrm{MF}}\max\{\bar R_{\max},1\}
\left[
\sum_{j=1}^K
(w_R^j+w_P^j)\alpha_{N_j}(E)
+
\sum_{j=1}^K
\frac{N_j}{N}\alpha_{N_j}(E)
\right], \label{eq.MFTerm}
\end{equation}
and
\begin{equation} \label{eq.HeterogeneityTerm}
  \epsilon_{\mathrm{heter}}
=
C_{\mathrm{heter}}
\left(
\sum_{t=0}^T\epsilon_t^R
+
\sum_{t=0}^{T-1}\epsilon_t^P
\right),  
\end{equation}
with
\[
\epsilon_t^R
:=
\left(
\frac{1}{N}
\sum_{i=1}^N
\sup_{\mathbf{s}\in S^N,\mathbf{a}\in A^N}
\left|
R_t^i(\mathbf{s},\mathbf{a})
-
\widehat R_t^i(\mathbf{s},\mathbf{a})
\right|^2
\right)^{1/2},
\]
and
\[
\epsilon_t^P
:=
\max_{i\in[N]}
\sup_{\mathbf{s}\in S^N,\mathbf{a}\in A^N}
\left\|
P_t^i(\cdot\mid\mathbf{s},\mathbf{a})
-
\widehat P_t^i(\cdot\mid\mathbf{s},\mathbf{a})
\right\|_{\mathrm{TV}},
\]
for constants $C_{\mathrm{MF}}=C_{\mathrm{MF}}(T,W_{\max},L_\pi,L_P,L_R,E)<\infty$ and $ C_{\mathrm{heter}}=C_{\mathrm{heter}}(T,W_{\max},L_\pi,L_P,L_R) <\infty$.
\end{theorem}

The two error terms have the following interpretations:
\begin{itemize}
    \item the \emph{mean-field approximation error} $\epsilon_{\mathrm{MF}}$ which measures the
    discrepancy between the homogenized $N$-player game
    $G_N^{\mathrm{MF}}$ and its $K$-population mean-field limit
    $G^{\mathrm{MF}}$. It involves the empirical measure estimate $\alpha_m(E)$ defined in \eqref{eq.alpham}.
    \item the \emph{heterogeneity error} $\epsilon_{\mathrm{heter}}$ which measures the discrepancy
    between the original game $G$ and the homogenized $N$-player game
    $G_N^{\mathrm{MF}}$.
\end{itemize}
The two error terms are analyzed separately in Theorems~\ref{thm:mf_approx} and~\ref{thm:perturbation}, whose combination yields Theorem~\ref{thm:homogenization}.


\subsection{Mean-field approximation for $N$-player games of mean-field type.}
We first bound the mean-field approximation error
$\epsilon_{\mathrm{MF}}$.
 
\begin{theorem}[Mean-field approximation]
\label{thm:mf_approx}
Let $G_N^{\mathrm{MF}}$ and $G^{\mathrm{MF}}$ be as in
Theorem~\ref{thm:homogenization}.
Let $\bar{\boldsymbol\pi}\in\Pi^K$ be a mean-field strategy profile
and let $\boldsymbol\pi\in\Pi^N$ denote its lift to
$G_N^{\mathrm{MF}}$, defined by
$$
\pi^i=\bar\pi^k,
\qquad i\in I_k.
$$
Under Assumptions~\ref{ass.boundedreward},  \ref{ass.weightedLipschitz}, \ref{ass.transition_lipschitz}, \ref{ass.reward_lipschitz}, and~\ref{ass.regularkernel}, there exists a constant
$
C_{\mathrm{MF}}=C_{\mathrm{MF}}(T,W_{\max},L_\pi,L_P,L_R,E)<\infty
$
such that
\[
\operatorname{AvgExpl}
(\boldsymbol\pi;G_N^{\mathrm{MF}})
\leq
\operatorname{AvgExpl}
(\bar{\boldsymbol\pi};G^{\mathrm{MF}})
+
\epsilon_{\mathrm{MF}},
\]
where
\[
\epsilon_{\mathrm{MF}}
=
C_{\mathrm{MF}}\max\{\bar R_{\max},1\}
\left[
\sum_{j=1}^K
(w_R^j+w_P^j)\alpha_{N_j}(E)
+
\sum_{j=1}^K
\frac{N_j}{N}\alpha_{N_j}(E)
\right].
\]
\end{theorem}
The proof is given in Section~\ref{proof.mf_approx}.

\begin{remark}[Scaling behavior of mean-field approximation error]
\label{rem:mf-error}
The mean-field approximation error in Theorem~\ref{thm:mf_approx}
is controlled by the empirical-measure modulus $\alpha_m(E)$:
$$
\epsilon_{\mathrm{MF}}
=
C_{\mathrm{MF}}\max\{\bar R_{\max},1\}
\left[
\sum_{j=1}^K
(w_R^j+w_P^j)\alpha_{N_j}(E)
+
\sum_{j=1}^K
\frac{N_j}{N}\alpha_{N_j}(E)
\right].
$$
Thus the dependence of the approximation error on the population
sizes is determined jointly by the sensitivity of the coefficients
to each population and by the rate of convergence of empirical
measures on the state-action space $E$.
For fixed $K$, since $\alpha_m(E)\to 0$ as ${m\to\infty}$ we have
$$\epsilon_{\mathrm{MF}}\to 0\qquad{\rm as}
\min_{1\leq j\leq K}N_j\to\infty$$ provided the weights
$w_R^j,w_P^j$ remain bounded.

A particularly important case is that in which the reward and
transition coefficients depend on the populations only through the
aggregate state-action distribution
$$
\bar L
=
\sum_{j=1}^K\frac{N_j}{N}L^j.
$$
If this dependence is respectively $w_R$- and $w_P$-Lipschitz, then
Assumption~\ref{ass.weightedLipschitz} holds with 
$$
w_R^j=w_R\frac{N_j}{N},
\qquad
w_P^j=w_P\frac{N_j}{N}.
$$
Consequently,
$$
\epsilon_{\mathrm{MF}}
\leq
C'
\sum_{j=1}^K
\frac{N_j}{N}\alpha_{N_j}(E),
$$
where $C'$ depends only on the constants appearing in
Theorem~\ref{thm:mf_approx}.
Using  $\alpha_m(E)\leq C_E m^{-1/d}$ 
then
$$
\epsilon_{\mathrm{MF}}
\leq
\frac{C'C_E}{N}
\sum_{j=1}^K N_j^{\frac{d-1}{d}}
\leq
C'C_E
\left(\frac{K}{N}\right)^{\frac{1}{d}}.
$$
The last inequality follows from the concavity of
$x\mapsto x^{1-1/d}$.
When $E$ is finite, sharper empirical-measure estimates recover the
usual square-root dependence on the population sizes.
\end{remark}

Note that in the classical mean-field approximation literature assuming
i.i.d. initial states, all of
$\max_{\widetilde{\pi}^i\in\Pi}
V^i(\pi^1,\dots,\widetilde{\pi}^i,\dots,\pi^N;G_N^{\mathrm{MF}})$
(resp. $V^i(\boldsymbol{\pi};G_N^{\mathrm{MF}})$),
$i\in I_k$, would be the same, and each of them is close to the
mean-field counterpart
$\max_{\widetilde{\pi}^k\in\Pi}
\bar{V}^k(\widetilde{\pi}^k,\bar{\boldsymbol{\pi}};G^{\mathrm{MF}})$
(resp.
$\bar{V}^k(\bar{\pi}^k,\bar{\boldsymbol{\pi}};G^{\mathrm{MF}})$).
As a result, one can bound individually each pair
\[
\left|
\max_{\widetilde{\pi}^k\in\Pi}
\bar{V}^k(\widetilde{\pi}^k,\bar{\boldsymbol{\pi}};G^{\mathrm{MF}})
-
\max_{\widetilde{\pi}^i\in\Pi}
V^i(\pi^1,\dots,\widetilde{\pi}^i,\dots,\pi^N;G_N^{\mathrm{MF}})
\right|
\]
and, respectively,
\[
\left|
\bar{V}^k(\bar{\pi}^k,\bar{\boldsymbol{\pi}};G^{\mathrm{MF}})
-
V^i(\boldsymbol{\pi};G_N^{\mathrm{MF}})
\right|.
\]
These methods are not directly applicable to our setup with heterogeneous deterministic initial states, since players within the same group may have different best-response values. In particular, directly bounding the aggregate best-response difference is infeasible because the summation and supremum cannot in general be exchanged due to the heterogeneity among players. To address this problem, we decompose the mean-field values according to the different initial states and compare each $V^i$ with the mean-field approximation corresponding to its own initial state $s_0^i$. This also leads to independent but not identically distributed comparison processes, for which we use Lemma~\ref{lem:empirical-BL}.

\begin{remark}
    The mean-field approximation  techniques are motivated by \cite{saldi2018markov}, where an asymptotic approximation result of the Nash equilibrium between mean-field games and competitive $N$-player games is established; non-asymptotic extensions have been studied by Yardim {\it et al.} \cite{yardim2024mean}. Our results generalize the analysis to account for heterogeneity both across populations and within the same population in terms of players' initial conditions. 
\end{remark}

\subsection{The contribution of population heterogeneity }
We now bound the heterogeneity error
$\epsilon_{\mathrm{heter}}$, which quantifies the discrepancy between the original heterogeneous $N$-player game $G$ and the homogenized $N$-player game of mean-field type $G_N^{\mathrm{MF}}$.

\begin{theorem}[Perturbation estimate]\label{thm:perturbation}
Let $G$ and $G_N^{\mathrm{MF}}$ be the $N$-player games introduced in
Theorem~\ref{thm:homogenization}.
Suppose that Assumptions~\ref{ass.boundedreward}, \ref{ass.weightedLipschitz},
\ref{ass.transition_lipschitz}, \ref{ass.reward_lipschitz},
and~\ref{ass.regularkernel} hold.
For $t\in\mathcal T$, define
\[
\epsilon_t^R
:=
\left(
\frac1N
\sum_{i=1}^N
\sup_{\mathbf{s}\in S^N,\mathbf{a}\in A^N}
\left|
R_t^i(\mathbf{s},\mathbf{a})
-
\widehat R_t^i(\mathbf{s},\mathbf{a})
\right|^2
\right)^{1/2},
\]
and, for $t=0,\ldots,T-1$, define
\[
\epsilon_t^P
:=
\max_{i\in[N]}
\sup_{\mathbf{s}\in S^N,\mathbf{a}\in A^N}
\left\|
P_t^i(\cdot\mid\mathbf{s},\mathbf{a})
-
\widehat P_t^i(\cdot\mid\mathbf{s},\mathbf{a})
\right\|_{TV}.
\]
Then there exists a constant $C_{\mathrm{heter}}=C_{\mathrm{heter}}(T,W_{\max},L_\pi,L_P,L_R)$
such that, for every admissible Markov strategy profile
$\bpi\in\Pi^N$,
\[
\operatorname{AvgExpl}(\boldsymbol{\pi};G)
-
\operatorname{AvgExpl}(\boldsymbol{\pi};G_N^{\mathrm{MF}})
\le
\epsilon_{\mathrm{heter}},
\]
where
\[
\epsilon_{\mathrm{heter}}
:=
C_{\mathrm{heter}}\left(
\sum_{t=0}^T\epsilon_t^R
+
\sum_{t=0}^{T-1}\epsilon_t^P
\right).
\]
\end{theorem}

\begin{remark}[Contribution of heterogeneity]
The reward and transition discrepancies are measured differently:
$\epsilon_t^R$ is an RMS average across players, whereas
$\epsilon_t^P$ is a worst-player discrepancy. Thus reward heterogeneity
is controlled on average, while transition heterogeneity is controlled
uniformly across players.
\end{remark}

\subsection{How accurate is the representative agent model?}\label{subsec:comparison2model}

Theorem \ref{thm:homogenization} shows that the  approximation error of a lifted mean-field strategy decomposes into two terms:

\begin{itemize}
    \item the mean-field approximation error $\epsilon_{\mathrm{MF}}$, and
    \item the impact  of heterogeneity  in each group: $\epsilon_{\mathrm{heter}}$.
\end{itemize} 
Intuitively, there is a  trade-off  between these two types of errors: as the number of sub-populations $K$ increases, the groups become smaller, so the mean-field approximation error generally increases, while the within-group heterogeneity decreases.
Consequently, it is expected that $\epsilon_{\mathrm{MF}}$ increases with $K$ while $\epsilon_{\mathrm{heter}}$ decreases.

We now illustrate this trade-off in the simplest  case   where all agents in the population are replaced by a single {\it representative agent}. This representative-agent paradigm is widely used in economics but has  also been  criticized as being inaccurate and misleading in many situations where heterogeneity and interactions play an important role \cite{kirman1992,gallegati1999}.
While these criticisms are largely methodological and qualitative,  they also have implications for the quantitative accuracy of the representative-agent approximation. Our results provide a way to quantify this approximation error.

Let us illustrate this phenomenon using    the example discussed in Section \ref{subsec:example}. 
The ``simplest'' heterogeneous case is where there are two different firm types,  consisting respectively of $\alpha N$ and $(1-\alpha)N$ firms, characterized by two different sets of coefficients $\{c_m^1\}_{m=0}^4$ and $\{c_m^2\}_{m=0}^4$ that differ only in one of the coefficients. Without loss of generality, we consider the case where $c^i_2$  is heterogeneous across agents, \ie, we have two possibilities:  $\alpha N$ firms with $\{c_0,c_1,c_2^1,c_3,c_4\}$ and  $(1-\alpha) N$ firms with $\{c_0,c_1,c_2^2,c_3,c_4\}$. 

We compare two models: first, a homogeneous representative agent model where all coefficients are replaced with a set of averaged coefficients $\{c_0,c_1,\bar{c}_2,c_3,c_4\}$ and, second, a heterogeneous model which incorporates this heterogeneity. Below we show that  while the first model entails a non-zero heterogeneity term $\epsilon_{\mathrm{heter}}$ compared to the second model, it exhibits a smaller mean-field approximation error $\epsilon_{\mathrm{MF}}$ due to averaging over a larger population.

\paragraph{Mean-field approximation error.}
In this example, the normalization $\operatorname{diam}(E)\leq1$
amounts to rescaling the natural metric by
$D_E:=S_{\max}+\sqrt{Q^2+H^2}$.

We assume that all cost coefficients above 
lie in the range $[0,C]$ for some constant $C>0$. In addition, assume $\sigma\in(0,1)$ and let $w:=
\frac{D_EQ\delta^{1/\sigma}}
{\sigma Q_0^{1+1/\sigma}}$. Then: 
\begin{itemize}
    \item For the first model, 
Assumption \ref{ass.weightedLipschitz} holds with $K=1$, and 
\[
\bar{R}_{\max}
    =
    \max\left\{
    (\delta/Q_0)^{1/\sigma}Q,\,
    C(Q+Q^2+H+2Q+S_{\max})
    \right\},
    \qquad
    w_R^1=w,
    \qquad
    w_P^1=0.
\]
Accordingly, $W_{\max}=w$.
\item Similarly, for the second model, Assumption \ref{ass.weightedLipschitz} also holds with the same $\bar{R}_{\max}$ and $W_{\max}$, while $w_R^1=\alpha w$,
    $w_R^2=(1-\alpha)w$, and
    $w_P^1=w_P^2=0$.
\end{itemize}  
Below we prove these two simple claims.
\begin{proof}{}
To see this, let us look at the first model as an example. Firstly, we have $K=1$ and that 
\begin{align*}
|\bar R_t^1(s_t,(q_t,h_t),L_t)|&=\left|p_tq_t^1-(c_0q_t^1+c_1(q_t^1)^2+c_2h_t^1+(c_2+c_3)\max\{q_t^1-s_t^1,0\}+c_4s_t^1)\right|\leq \bar{R}_{\max},
\end{align*} where we make use of the fact that the  $p_t\leq (\delta/Q_0)^{1/\sigma}$ and hence $0\leq p_tq_t^1\leq (\delta/Q_0)^{1/\sigma}Q$, and 
\[
0\leq c_0q_t^1+c_1(q_t^1)^2+c_2h_t^1+(c_2+c_3)\max\{q_t^1-s_t^1,0\}+c_4s_t^1\leq CQ+CQ^2+CH+2CQ+CS_{\max}.
\]
In addition, we have 
\[
\left|\bar{R}_t(s_t^1,a_t^1,L_t^1)-\bar{R}_t(s_t^1,a_t^1,\widetilde{L}_t^1)\right|=|p_t-\widetilde{p}_t|q_t^1,
\]
where $p_t=\left(\frac{\delta}{Q_0+\sum_{q}q\sum_{s,h}L_t^1(s,q,h)}\right)^{1/\sigma}$, and $\widetilde{p}_t=\left(\frac{\delta}{Q_0+\sum_{q}q\sum_{s,h}\widetilde{L}_t^1(s,q,h)}\right)^{1/\sigma}$, and hence 
\begin{align*}
\left|\bar{R}_t(s_t^1,a_t^1,L_t^1)-\bar{R}_t(s_t^1,a_t^1,\widetilde{L}_t^1)\right|&\leq \frac{Q\delta^{1/\sigma}}
{\sigma Q_0^{1+1/\sigma}}\left|\sum_{q}q\sum_{s,h}L_t^1(s,q,h)-\sum_{q}q\sum_{s,h}\widetilde{L}_t^1(s,q,h)\right|\\
&\leq \frac{D_EQ\delta^{1/\sigma}}
{\sigma Q_0^{1+1/\sigma}}W_1(L_t^1,\widetilde{L}_t^1),
\end{align*}
where we make use of the mean-value theorem and the fact that $|f'(x)|\leq \frac{\delta^{1/\sigma}}{\sigma Q_0^{1+1/\sigma}}$ for $x\geq 0$ where $f(x):=(\delta/(Q_0+x))^{1/\sigma}$, as well as $\left|
\sum_q q\sum_{s,h}L_t^1(s,q,h)
-
\sum_q q\sum_{s,h}\widetilde L_t^1(s,q,h)
\right|
\leq
D_E W_1(L_t^1,\widetilde L_t^1)$, which follows from the normalized metric above.

The fact that $w_P^1=0$ is obvious as the transition model is mean-field independent here. 
The derivation for the second model is nearly identical and is hence omitted. 
\end{proof}

Now by plugging these quantities in Theorem \ref{thm:homogenization}, 
we can compare the mean-field approximation errors for the two models.
Since $E$ is finite, Lemma~\ref{lem:empirical-BL} gives
\[
\alpha_m(E)=\sqrt{\frac{|E|}{m}},
\qquad
|E|=(S_{\max}+1)(Q+1)(H+1).
\]
Hence, for the first model,
\[
\epsilon_{\mathrm{MF}}^1
=
\widetilde C\frac{w+1}{\sqrt N},
\qquad
\text{where}\quad
\widetilde C
:=
C_{\mathrm{MF}}\max\{\bar R_{\max},1\}\sqrt{|E|},
\]
with the same constant $C_{\mathrm{MF}}$ as in Theorem~\ref{thm:homogenization}.

For the second model, with
$N_1=\alpha N$ and $N_2=(1-\alpha)N$, we obtain
\[
\begin{aligned}
\epsilon_{\mathrm{MF}}^2
&=
\frac{\widetilde C}{\sqrt N}
\left[
w\big(\sqrt{\alpha}+\sqrt{1-\alpha}\big)
+
\sqrt{\alpha}+\sqrt{1-\alpha}
\right]
=
\widetilde C\frac{w+1}{\sqrt N}
\big(\sqrt{\alpha}+\sqrt{1-\alpha}\big).
\end{aligned}
\]
Therefore,
\[
\frac{\epsilon_{\mathrm{MF}}^2}
{\epsilon_{\mathrm{MF}}^1}
=
\sqrt{\alpha}+\sqrt{1-\alpha}.
\]
This indicates that model 2, which better identifies the heterogeneity among the two groups of firms, leads to a certified mean-field approximation bound that is larger by
$\sqrt{\alpha}+\sqrt{1-\alpha}-1$ in relative terms compared to the first model that cannot distinguish the two groups of firms.

\paragraph{Heterogeneity error.} As the second model correctly identifies the two groups of firms,  the heterogeneity error $\epsilon_{\mathrm{heter}}^2$   as defined in Theorem \ref{thm:homogenization} is zero for this model. In contrast, the heterogeneity error $\epsilon_{\mathrm{heter}}^1$ for the first model depends on the dispersion of the coefficients, and is given by $C_{\mathrm{heter}}\sum_{t=0}^T\epsilon_t^R$, where $\alpha^1=\alpha$, $\alpha^2=1-\alpha$, and  
\begin{align*}
\epsilon_t^R&=\sqrt{\sum_{i=1,2}\alpha^i\max_{s_t^i\in\{0,\dots,S_{\max}\},q_t^i\in\{0,\dots,Q\},h_t^i\in\{0,\dots,H\}}|(c_2^i-\bar{c}_2)h_t^i+(c_2^i-\bar{c}_2)\max\{q_t^i-s_t^i,0\}|^2}\\
&=\max_{s\in\{0,\dots,S_{\max}\},q\in\{0,\dots,Q\},h\in\{0,\dots,H\}}|h+\max\{q-s,0\}|\sqrt{\alpha|c_2^1-\bar c_2|^2
+
(1-\alpha)|c_2^2-\bar c_2|^2}\\
&=(H+Q)\sqrt{\alpha|c_2^1-\bar{c}_2|^2+(1-\alpha)|c_2^2-\bar{c}_2|^2},
\end{align*}
and we make use of the fact that $\epsilon_t^P=0$ as we make no approximation in the transition model, which is already homogeneous across the firms by definition. As a result, we have $\epsilon_{\mathrm{heter}}^1=C_{\mathrm{heter}}(H+Q)(T+1)\sqrt{\alpha|c_2^1-\bar{c}_2|^2+(1-\alpha)|c_2^2-\bar{c}_2|^2}$, which, among different values of $\bar{c}_2$, is minimized at $\bar{c}_2=\alpha c_2^1+(1-\alpha)c_2^2$. The resulting heterogeneity error of the first model is then $\epsilon_{\mathrm{heter}}^1=C_{\mathrm{heter}}(H+Q)(T+1)\sqrt{\alpha(1-\alpha)}|c_2^1-c_2^2|$ (as opposed to the zero heterogeneity error of the second model). 

\paragraph{Trade-off between two sources of errors.}
When the two sources of errors are combined, it is easy to see that the trade-off between these two errors comes from the trade-off between the number of firms $N$ and the heterogeneity among the firms (the difference between $c_2^1$ and $c_2^2$).

For small $N$, the reduction in mean-field approximation error obtained by using a single larger population may outweigh the heterogeneity error. For sufficiently large $N$, however, the mean-field errors become small and the nonzero heterogeneity error of the single-population model dominates; the two-population model is then more accurate. For a fixed $N$, greater heterogeneity favors the two-population model. 

Note that $\frac{\epsilon_{\mathrm{MF}}^2}{\epsilon_{\mathrm{MF}}^1}
=
\sqrt{\alpha}+\sqrt{1-\alpha}<2$. Therefore, if $\epsilon_{\mathrm{MF}}^1
<
\epsilon_{\mathrm{heter}}^1$, then $\epsilon_{\mathrm{heter}}^1+\epsilon_{\mathrm{MF}}^1
>
2\epsilon_{\mathrm{MF}}^1
>
\epsilon_{\mathrm{MF}}^2$, which means that the single-population MFG has a larger certified approximation-error bound.

Using
\[
\epsilon_{\mathrm{MF}}^1
=
\frac{\widetilde C(w+1)}{\sqrt N}, \qquad 
\epsilon_{\mathrm{heter}}^1
=
C_{\mathrm{heter}}(H+Q)(T+1)\sqrt{\alpha(1-\alpha)}
|c_2^1-c_2^2|,
\]
the above sufficient condition holds whenever
\[
N>N^*:=
\frac{\bar K}
{\alpha(1-\alpha)|c_2^1-c_2^2|^2},
\]
where
\[
\bar K
:=
\left(
\frac{\widetilde C(w+1)}
{C_{\mathrm{heter}}(H+Q)(T+1)}
\right)^2.
\]
Therefore, if the number of firms $N$ exceeds the threshold $N^*$, 
the homogeneous representative-firm model has a larger certified
approximation-error bound than the two-population mean-field model.

The variability within the population, which influences this bound, depends on two key factors: the proportion of each type of firm, $\alpha$ and $1-\alpha$, and the difference between the parameter sets associated with each firm type, $|c_2^1 - c_2^2|$. The term $\alpha(1-\alpha)$ captures the composition of the population, attaining its maximum when $\alpha = 0.5$, where the two firm types are most evenly mixed. As $\alpha$ approaches either 0 or 1, the population becomes more homogeneous, reducing the denominator and thereby increasing the sufficient threshold $N^*$. Similarly, the term $|c_2^1 - c_2^2|$ measures the distinction between the two groups; when firms are nearly identical, this difference is small, leading to a larger $N^*$, indicating that a single-population MFG may remain a reasonable approximation for a wider range of $N$. Conversely, when the two types of firms are more distinct, the denominator increases, lowering $N^*$, meaning that even for relatively small $N$, the two-population MFG may provide a  better representation. Therefore, a more homogeneous population -- either due to an imbalanced composition ($\alpha$ near 0 or 1) or a small difference in firm characteristics -- leads to  a larger sufficient threshold $N^*$, making the single-population MFG favorable for a broader range of system sizes. In contrast, greater heterogeneity, characterized by a balanced mix of firms and substantial differences between their parameters, lowers $N^*$, favoring a two-population MFG even for moderate $N$. This relationship underscores the role of population diversity in determining the appropriate modeling approach in multi-agent games.

\section{Partitioning the population of players}\label{sec:partition}

There are many ways to partition the population of players into groups which may then be approximated as
homogeneous, and some partitions may lead to substantially better approximations than others. We now use
the quantitative homogenization estimates of Section~\ref{sec:homogenization_error} to compare different
multi-population mean-field approximations and to optimize the choice of the partition.

Throughout this section we work with the compact state-action space
$E:= S\times A\subset\mathbb R^{d}$ and assume $d\geq 3$.
We retain Assumptions of Section~\ref{sec:homogenization_error},
uniformly in the heterogeneity parameter introduced below. Recall from Lemma \ref{lem:empirical-BL} that 
\[
\alpha_m(E)\leq C_E m^{-\frac{1}{d}}
\]
is an upper  bound for the expected Wasserstein-1 error of an empirical measure based on
$m$ independent $E$-valued samples. 
For $d=1$ or $2$, the sharper bounds in Lemma~\ref{lem:empirical-BL} may be used instead, with the corresponding modifications to the rates below.

\paragraph{A parametrized heterogeneous model.}
We consider an $N$-player game
\[
G=( S, A,\{P_t^i,R_t^i\}_{t\in\mathcal T,i\in[N]},{\bf s}_0)
\]
whose reward functions depend on player-specific characteristics
$
\theta^i\in B_{d_\theta}(0,D):= \{\theta\in\mathbb R^{d_\theta}:\|\theta\|_2\le D\}\subset\mathbb R^{d_\theta}.
$ We regard the characteristics $\theta^i$ as given inputs to the partitioning problem, assumed to be available either through direct observation or prior estimation.
More precisely,
$
R_t^i({\bf s},{\bf a})
=
\bar R_t\bigl(s^i,a^i,L^N({\bf s},{\bf a}),\theta^i\bigr),
$
where
\[
L^N({\bf s},{\bf a})
=
\frac1N\sum_{j=1}^N\delta_{(s^j,a^j)}\in\mathcal P(E).
\]
For simplicity, the transition dynamics are assumed to be homogeneous across players:
\[
P_t^i(\,\cdot\,|{\bf s},{\bf a})
=
\bar P_t\bigl(\,\cdot\,|s^i,a^i,L^N({\bf s},{\bf a})\bigr).
\]
Thus, in this section, heterogeneity is carried by the reward parameters $\theta^i$. The same construction may be combined with Theorem~\ref{thm:perturbation} when transition kernels are also heterogeneous, under an analogous Lipschitz condition with respect to the heterogeneity parameter.

We impose the following parameter regularity.

\begin{assumption}[Parametric Lipschitz continuity]\label{parametrized_setup_assumption}
There exist constants $w_R,w_P,w_D\geq0$ such that, for all $t\in\mathcal T$,
$(s,a)\in E$, $L,\widetilde L\in\mathcal P(E)$ and
$\theta,\widetilde\theta\in B_{d_\theta}(0,D)$,
\[
\left|
\bar R_t(s,a,L,\theta)-\bar R_t(s,a,\widetilde L,\theta)
\right|
\leq
w_R\,W_1(L,\widetilde L),
\]
\[
\left\|
\bar P_t(\,\cdot\,|s,a,L)
-
\bar P_t(\,\cdot\,|s,a,\widetilde L)
\right\|_{\mathrm{TV}}
\leq
w_P\,W_1(L,\widetilde L),
\]
and
\[
\left|
\bar R_t(s,a,L,\theta)-\bar R_t(s,a,L,\widetilde\theta)
\right|
\leq
w_D\|\theta-\widetilde\theta\|_2.
\]
The reward functions are uniformly bounded by $\bar R_{\max}<\infty$.
\end{assumption}
Let
$
W_{\max}:=\max\{w_R,w_P\}.
$
\paragraph{Mean-field approximations associated with a partition.}
Let
\[
[N]=\bigsqcup_{k=1}^K  I_k,
\qquad
N_k:=|  I_k|,
\]
be a partition into nonempty groups. We homogenize group $k$ by replacing the individual
parameters $\theta^i$, $i\in  I_k$, with the centroid
\[
\widehat\theta^k
:=
\frac1{N_k}\sum_{i\in  I_k}\theta^i.
\]
The corresponding $K$-population MFG is denoted 
$
G^{\mathrm{MF}}(\{  I_k\}_{k\in[K]}).
$
Its transition and reward kernels are
\[
\bar P_t^k(\,\cdot\,|s,a,L^1,\ldots,L^K)
:=
\bar P_t(\,\cdot\,|s,a,L),\quad
\bar R_t^k(s,a,L^1,\ldots,L^K)
:=
\bar R_t(s,a,L,\widehat\theta^k),
\]
where
\[
L:=\sum_{k=1}^K\frac{N_k}{N}L^k.
\]
The initial distribution of population $k$ is the empirical measure
\[
\mu_0^{k,N}
:=
\frac1{N_k}\sum_{i\in  I_k}\delta_{s_0^i}.
\]

Since
\[
W_1\left(
\sum_{k=1}^K\frac{N_k}{N}L^k,
\sum_{k=1}^K\frac{N_k}{N}\widetilde L^k
\right)
\leq
\sum_{k=1}^K\frac{N_k}{N}W_1(L^k,\widetilde L^k),
\]
Assumption~\ref{parametrized_setup_assumption} implies that the weighted Lipschitz constants in
Theorem~\ref{thm:mf_approx} may be chosen as
\[
w_R^k=w_R\frac{N_k}{N},
\qquad
w_P^k=w_P\frac{N_k}{N},
\qquad k\in[K].
\]

The next result is the specialization of Theorems~\ref{thm:mf_approx} and~\ref{thm:perturbation} to this
parametric setting.
\begin{corollary}[Partition-dependent homogenization bound]
\label{cor:partition-bound}
Let $G^{\mathrm{MF}}(\{I_k\}_{k\in[K]})$ be the $K$-population mean-field game associated with the
partition
$$
[N]=\bigsqcup_{k=1}^K I_k,
\qquad
N_k:=|I_k|,
$$
and
$
\bar{\boldsymbol\pi}
=
(\bar\pi^1,\ldots,\bar\pi^K)\in\Pi^K
$
 a mean-field strategy profile satisfying
$$
\operatorname{AvgExpl}
\left(
\bar{\boldsymbol\pi};
G^{\mathrm{MF}}(\{I_k\}_{k\in[K]})
\right)
\leq \epsilon,
$$
and define its lift $\boldsymbol\pi\in\Pi^N$ by
$\pi^i=\bar\pi^k, i\in I_k.$
Under Assumptions \ref{ass.boundedreward}, \ref{ass.weightedLipschitz}, \ref{ass.transition_lipschitz}, \ref{ass.reward_lipschitz}, \ref{ass.regularkernel},   and ~\ref{parametrized_setup_assumption} we have
$$
\operatorname{AvgExpl}(\boldsymbol\pi;G)
\leq
\epsilon
+
\epsilon_{\mathrm{MF}}(\{I_k\}_{k\in[K]})
+
\epsilon_{\mathrm{heter}}(\{I_k\}_{k\in[K]}),
$$
where
\begin{equation}
 \epsilon_{\mathrm{MF}}(\{I_k\}_{k\in[K]})
:=
C_0
\sum_{k=1}^K
\frac{N_k}{N}\,
\alpha_{N_k}(E),   \label{epsilon_MF_partition}
\end{equation}
where $
C_0
:=
C_{\mathrm{MF}}
\max\{\bar R_{\max},1\}
(w_R+w_P+1),
$
and
\begin{equation}
\epsilon_{\mathrm{heter}}(\{I_k\}_{k\in[K]})
:=
C_D
\left(
\frac1N
\sum_{k=1}^K
\sum_{i\in I_k}
\|\theta^i-\widehat\theta^k\|_2^2
\right)^{1/2}, 
 \label{epsilon_heter_partition}   
\end{equation}
where $C_D
:=
C_{\mathrm{heter}}(T+1)w_D$. 
In particular, for $E\subset\mathbb R^{d}$ with $d\geq 3$,
$\alpha_m(E)\leq C_E m^{-1/d},$
so
\begin{equation}
\epsilon_{\mathrm{MF}}(\{I_k\}_{k\in[K]})
\leq
\frac{C_0C_E}{N}
\sum_{k=1}^K
N_k^{1-\frac{1}{d}}.    \label{epsilon_approx_partition}
\end{equation}
\end{corollary}
\textit{Proof.}
For the aggregate mean field
$L=\sum_k(N_k/N)L^k$, the weighted Lipschitz coefficients appearing in
Theorem~\ref{thm:mf_approx} are
\[
w_R^k=w_R\frac{N_k}{N},
\qquad
w_P^k=w_P\frac{N_k}{N}.
\]
Substituting these coefficients in Theorem~\ref{thm:mf_approx} gives
\[
\epsilon_{\mathrm{MF}}
\leq
C_{\mathrm{MF}}\max\{\bar R_{\max},1\}
\left[
\sum_{k=1}^K(w_R+w_P)\frac{N_k}{N}\alpha_{N_k}(E)
+
\sum_{k=1}^K\frac{N_k}{N}\alpha_{N_k}(E)
\right],
\]
which yields \eqref{epsilon_MF_partition}.
For the perturbation step there is no transition heterogeneity. Moreover,
for $i\in  I_k$,
\[
\sup_{{\bf s}\in S^N,\ {\bf a}\in A^N}
\left|
R_t^i({\bf s},{\bf a})-\widehat R_t^i({\bf s},{\bf a})
\right|
\leq
w_D\|\theta^i-\widehat\theta^k\|_2.
\]
Hence the reward-heterogeneity term satisfies
\[
\epsilon_t^R
\leq
w_D
\left(
\frac1N\sum_{k=1}^K\sum_{i\in  I_k}
\|\theta^i-\widehat\theta^k\|_2^2
\right)^{1/2}.
\]
Since there is no transition heterogeneity, Theorem~\ref{thm:perturbation} therefore gives 
\eqref{epsilon_heter_partition}. The last statement follows by using
$\alpha_m(E)\leq C_E m^{-1/d}$.

\begin{remark}[Augmenting the state by the heterogeneity parameter]
An alternative is to regard $\theta^i$ as part of the state and work on an augmented space
$ S\times B_{d_\theta}(0,D)$. In the compact-space formulation this does not merely
increase a finite cardinality: it increases the metric dimension of the empirical-measure problem.
Consequently, the generic bound $\alpha_m(E)$ may deteriorate. For example, if the effective
state-action dimension increases from $d$ to $d+d_\theta$, the generic Euclidean rate changes
from $m^{-1/d}$ to $m^{-1/(d+d_\theta)}$. Population partitioning therefore provides a
different approximation trade-off: it keeps the original state-action dimension and pays instead the
explicit heterogeneity term \eqref{epsilon_heter_partition}.
\end{remark}

We now study the \emph{bound-optimal} partitioning problem
\begin{equation}\label{eq:total_error}
\min_{\{  I_k\}}
\left\{
\epsilon_{\mathrm{MF}}(\{  I_k\})
+
\epsilon_{\mathrm{heter}}(\{  I_k\})
\right\}.
\end{equation}
Thus ``optimal'' below always refers to the certified upper bound on average exploitability supplied by Corollary \ref{cor:partition-bound}, rather than to the  exact average exploitability of the lifted profile.

\subsection{Optimal partitioning of a large population}

For a fixed number of groups $K$, the compact-space mean-field term satisfies a useful
partition-independent upper bound. Using 
$\alpha_m(E)\leq C_Em^{-1/d}$ and the  concavity of $x\mapsto x^{1-1/d}$,
\begin{align}
\epsilon_{\mathrm{MF}}(\{  I_k\}_{k\in[K]})
&\leq
\frac{C_0C_E}{N}\sum_{k=1}^K N_k^{1-\frac{1}{d}}
\leq
C_0C_E
\left(\frac{K}{N}\right)^{\frac{1}{d}}.
\label{eq:epsilon_approx_to_0}
\end{align}
In particular, for fixed $K$,
$\epsilon_{\mathrm{MF}}\to0$ as $N\to\infty$.
More generally, the bound vanishes whenever $K=K_N=o(N)$.

For fixed $K$, minimizing the heterogeneity term
\eqref{epsilon_heter_partition} is equivalent to the $K$-means clustering objective:
\begin{equation}\label{epsilon_heter_partition_group_k}
\begin{array}{ll}
\displaystyle
\text{\rm minimize}_{\{  I_k\}_{k\in[K]}:\,
\bigsqcup_{k=1}^K  I_k=[N]}
&
\displaystyle
D_{N,K}:=
\frac1N
\sum_{k=1}^K\sum_{i\in  I_k}
\|\theta^i-\widehat\theta^k\|_2^2,
\\[2mm]
&
\displaystyle
\widehat\theta^k=
\frac1{|  I_k|}
\sum_{i\in  I_k}\theta^i.
\end{array}
\end{equation}
$D_{N,K}$ is precisely the classical K-means objective.
 
\paragraph{Case study: uniformly distributed agent features.}
Assume now that $d_\theta=1$ and that
$\theta^1,\ldots,\theta^N$ are i.i.d.\ $\operatorname{Uniform}[0,1]$.
For fixed $K$, the population $K$-means problem has the limiting centroids
\[
\frac{2k-1}{2K},
\qquad k=1,\ldots,K,
\]
and, by the consistency of optimal quantization/K-means
\cite[Section~2]{pollard1982quantization},
\begin{equation}\label{eq:uniform_kmeans_distortion}
D_{N,K}\longrightarrow \frac{1}{12K^2}
\qquad\text{almost surely as }N\to\infty.
\end{equation}
Consequently,
\begin{equation}\label{eq:epsilon_heter_K}
\epsilon_{\mathrm{heter}}
\longrightarrow
\frac{C_D}{\sqrt{12}}\frac1K
\qquad\text{almost surely}.
\end{equation}
Thus, in this example, splitting a large population into $K$ optimally chosen groups reduces the
heterogeneity contribution by a factor asymptotic to $1/K$, while the mean-field contribution grows
with $K$.

The compact-space estimate also makes the dimensional trade-off explicit. Suppose, for the purpose
of a rate calculation, that $K=K_N$ grows slowly enough that the $K$-means groups remain
approximately balanced, $N_k\asymp N/K$, and that the quantization distortion retains the order
$D_{N,K}\asymp K^{-2}$. Then
\[
\epsilon_{\mathrm{MF}}
\lesssim
\left(\frac{K}{N}\right)^{\frac{1}{d}},
\qquad
\epsilon_{\mathrm{heter}}
\asymp
\frac1K.
\]
Balancing these two rates suggests
\begin{equation}\label{eq:optimal_K_scaling}
K
\asymp
N^{1/(d+1)},
\end{equation}
and suggests a total bound  of order
\[
N^{-\frac{1}{d+1}}.
\]
This calculation should be understood as the scaling suggested by the two error bounds; a precise
growing-$K$ statement would additionally require uniform quantization estimates.

\subsection{Bound-optimal partitioning for finite populations}

For finite $N$, minimizing the heterogeneity term alone need not minimize the full certified error
bound, because the mean-field term penalizes small groups. We now give an exact mixed-integer
second-order cone formulation of the bound-optimal partitioning problem.

Fix an upper bound $K_{\max}\leq N$ on the number of groups. Empty group slots are allowed, so the
number of active groups is itself selected by the optimization. Define
\[
\psi_E(0):=0,
\qquad
\psi_E(m):=m\,\alpha_m(E),
\quad m=1,\ldots,N.
\]
For the Euclidean bound from Lemma~\ref{lem:empirical-BL}, 
$
\psi_E(m)=C_E m^{1-\frac{1}{d}}=C_E m^{\frac{d-1}{d}}.
$
\begin{theorem}[Mixed-integer formulation of the bound-optimal partition]
\label{thm:MICP_optimal_agent_partitioning}
Let $K_{\max}\leq N$. The minimization of
$
\epsilon_{\mathrm{MF}}+\epsilon_{\mathrm{heter}}
$
over all partitions of $[N]$ into at most $K_{\max}$ nonempty groups is equivalent to the following
mixed-integer second-order cone program. The variables are
\[
w_{ik}\in\{0,1\},
\quad
z_{km}\in\{0,1\},
\quad
d_{ik}\geq0,
\quad
c_k\in\mathbb R^{d_\theta},
\quad
x\geq0,
\]
for $i\in[N]$, $k\in[K_{\max}]$, and $m\in\{0,\ldots,N\}$:
\begin{align}
\text{\rm minimize}\quad
&
\frac{C_D}{\sqrt N}\,x
+
\frac{C_0}{N}
\sum_{k=1}^{K_{\max}}\sum_{m=0}^N
\psi_E(m)z_{km}
\label{formulation:MICP}
\\
\text{\rm subject to}\quad
&
\sum_{k=1}^{K_{\max}}w_{ik}=1,
\qquad i\in[N],
\notag\\
&
\sum_{m=0}^N z_{km}=1,\quad \sum_{i=1}^Nw_{ik}
=
\sum_{m=0}^N m z_{km},
\quad k\in[K_{\max}],
\notag\\
&
d_{ik}
\geq
\|\theta^i-c_k\|_2-2D(1-w_{ik}),
\qquad i\in[N],\ k\in[K_{\max}],
\notag\\
&
\|c_k\|_2\leq D,
\qquad k\in[K_{\max}],
\notag\\
&
\left(
\sum_{k=1}^{K_{\max}}\sum_{i=1}^N d_{ik}^2
\right)^{1/2}
\leq x.
\notag
\end{align}
\end{theorem}

\begin{remark}
The variables $z_{km}$ encode the cardinality of group $k$. Their use is convenient because the
compact-space mean-field penalty depends nonlinearly on the group size through
$m\alpha_m(E)$, while the values $\psi_E(m)$ are fixed coefficients for
$m=0,\ldots,N$. Hence the mean-field part of the objective is linear in the size-selection
variables. The remaining nonlinear constraints are second-order conic.
\end{remark}
\textit{Proof.}
Given a feasible assignment $\{w_{ik}\}$, define
\[
  I_k:=\{i\in[N]:w_{ik}=1\},
\qquad
N_k:=|  I_k|.
\]
The size-selection constraints force $z_{k,N_k}=1$ and $z_{km}=0$ for $m\neq N_k$. Therefore,
\[
\frac{C_0}{N}
\sum_{k=1}^{K_{\max}}\sum_{m=0}^N\psi_E(m)z_{km}
=
C_0
\sum_{k:N_k>0}\frac{N_k}{N}\alpha_{N_k}(E)
=
\epsilon_{\mathrm{MF}}.
\]
For a fixed assignment, minimization over $c_k$ of
\[
\sum_{i\in  I_k}\|\theta^i-c_k\|_2^2
\]
is attained at the centroid
\[
c_k=\frac1{N_k}\sum_{i\in  I_k}\theta^i
\]
whenever $N_k>0$. Since the ball $B_{d_\theta}(0,D)$ is convex and
$\theta^i\in B_{d_\theta}(0,D)$, this centroid satisfies $\|c_k\|_2\leq D$.

The big-$M$ constraint is exact because
\[
\|\theta^i-c_k\|_2\leq2D.
\]
Thus, for any fixed assignment, we may choose
\[
d_{ik}
=
\begin{cases}
\|\theta^i-c_k\|_2,& w_{ik}=1,\\
0,&w_{ik}=0.
\end{cases}
\]
The second-order cone constraint may then be taken to be tight, giving
\[
x
=
\left(
\sum_{k:N_k>0}\sum_{i\in  I_k}
\|\theta^i-\widehat\theta^k\|_2^2
\right)^{1/2}.
\]
Consequently,
\[
\frac{C_D}{\sqrt N}x
=
\epsilon_{\mathrm{heter}}.
\]
The objective in \eqref{formulation:MICP} is therefore exactly the partition-dependent upper bound
in Corollary~\ref{cor:partition-bound}. Conversely, every partition into at most
$K_{\max}$ groups defines feasible variables with the same objective value by taking
$w_{ik}={\bf1}_{\{i\in  I_k\}}$, choosing the corresponding size variable
$z_{k,N_k}=1$, and taking $c_k$ equal to the centroid. This proves the equivalence.

\begin{remark}
The formulation in Theorem~\ref{thm:MICP_optimal_agent_partitioning} optimizes the certified
homogenization bound, not the exact exploitability of the MFG equilibrium in the original game.
This distinction is useful in interpreting the resulting partition. The formulation also accommodates
sharper empirical-measure estimates: if Theorem~\ref{thm:mf_approx} is improved by replacing
$\alpha_m(E)$ with another known rate, only the precomputed coefficients
$\psi_E(m)=m\alpha_m(E)$ need to be changed.
\end{remark}

\section{Conclusion}

We developed a quantitative framework for approximating heterogeneous multi-player games by multi-population mean-field games. The approximation proceeds through an intermediate finite-player game of mean-field type, leading to a non-asymptotic bound on average exploitability that separates the finite-population mean-field approximation error from the error due to within-group heterogeneity. The framework allows heterogeneity in rewards, transition dynamics, and deterministic initial states, and applies to arbitrary partitions of the player population.

These bounds also provide a quantitative criterion for choosing the population partition, balancing the reduction of within-group heterogeneity against the increase in finite-population mean-field error caused by smaller groups. In the parametrized setting considered in Section~5, this leads to a finite-population mixed-integer second-order cone formulation, while in the large-population regime the heterogeneity component is closely related to $K$-means clustering. Natural directions for future work include learning the relevant heterogeneity parameters and partitions from data, incorporating richer forms of heterogeneous dynamics into the partitioning problem, and developing computational methods for large-scale instances.

\paragraph{Acknowledgment.}
A.H.'s research was funded through a Hooke Fellowship from the Oxford Mathematical Institute.  Part of this research was completed during our visit to the Institute for Mathematical and Statistical Innovation (IMSI), supported by the National Science Foundation (Grant No. DMS-1929348). R.C. thanks Ioannis Karatzas and Columbia University's Department of Mathematics for their kind hospitality, which enabled us to finish this project.

\newpage
\section{Technical proofs }
\subsection{Proof of Lemma~\ref{lem:lipschitz_policy_approximation}}
Fix $k\in[K]$ and a mean-field strategy profile $\bar{\bm\pi}$. The population flow $\left(L_t^{1,\bar{\bm\pi}},\ldots,L_t^{K,\bar{\bm\pi}}\right)_{t\in\mathcal T}$ is fixed in the unilateral optimization problem of a representative player in population $k$. Define the optimal value functions and the optimal state-action value functions recursively by $v_{T+1}(s):=0$, and for $t=T,\dots,0$,
\[
Q_t(s,a):=\bar R_t^k
\left(s,a,L_t^{1,\bar{\bm\pi}},\dots,L_t^{K,\bar{\bm\pi}}\right)+\mathbf 1_{\{t<T\}}\int_S
v_{t+1}(s')
\bar P_t^k\left(ds'\mid s,a, L_t^{1,\bar{\bm\pi}},
\dots,L_t^{K,\bar{\bm\pi}}\right),
\]
and
\[
v_t(s):=\sup_{a\in  A}Q_t(s,a).
\]
By Assumptions~\ref{ass.transition_lipschitz} and \ref{ass.reward_lipschitz} and backward induction, 
\[
|Q_t(s,a)-Q_t(s',a')|\le C_t\bigl(|s-s'|+|a-a'|\bigr)
\]
and 
\[
 |v_t(s)-v_t(s')|\le C_t |s-s'|,
\]
where the constants $C_t$ are independent of $k$ and $\bar{\bm\pi}$ and may be chosen recursively as  $C_{T+1}:=0$ and $C_t:=L_R+L_P C_{t+1}$ for $t=T,\dots,0$. 

Fix $\delta>0$. Since $S$ is a compact subset of a finite-dimensional Euclidean space, let $\{s_i\}_{i=1}^M$ be a maximal $\delta$-separated subset of $S$. By maximality, $\{s_i\}_{i=1}^M$ is also a $\delta$-net of $S$.  For each $t\in\mathcal T$ and $i\in[M]$, choose 
\[ 
a_{t,i}\in\arg\max_{a\in  A}Q_t(s_i,a),
\]
which exists by the compactness of $A$ and the continuity of $Q_t(s_i,\cdot)$.
We interpolate these locally optimal actions using the weights 
\[ 
\psi_i(s) := \left(1-\frac{|s-s_i|}{2\delta}\right)_+, \qquad \varphi_i(s) := \frac{\psi_i(s)} {\sum_{j=1}^M\psi_j(s)}. 
\]
Since $\{s_i\}$ is a $\delta$-net, the denominator is bounded below by $1/2$. Moreover, 
\[ 
\varphi_i(s)>0 \quad\Longrightarrow\quad |s-s_i|<2\delta. 
\]
Since the points $\{s_i\}$ are $\delta$-separated and $ S$ is contained in a finite-dimensional Euclidean space, 
at every $s\in S$, the number of nonzero $\varphi_i(s)$ is uniformly bounded by some constant $C_S<\infty$, independent of $s$ and $\delta$. 
Together with the $O(\delta^{-1})$ Lipschitz continuity of $\psi_i$, this yields 
\[ 
\sum_{i=1}^M |\varphi_i(s)-\varphi_i(s')| \le \frac{C_S}{\delta}|s-s'|, \qquad s,s'\in S, 
\]
for a constant $C_S$ depending only on $S$. Define the policy 
\[
\pi_t^\delta(\cdot\mid s) := \sum_{i=1}^M \varphi_i(s)\delta_{a_{t,i}}.\]  
Since 
$\operatorname{diam}(A)\le \operatorname{diam}(E)\le 1$,
\[ 
W_1\bigl( \pi_t^\delta(\cdot\mid s), \pi_t^\delta(\cdot\mid s') \bigr) \le \sum_{i=1}^M |\varphi_i(s)-\varphi_i(s')| \le \frac{C_S}{\delta}|s-s'|, 
\] so that $\pi^\delta\in\Pi_{C_S/\delta}.$ It remains to bound the loss of optimality. If $\varphi_i(s)>0$, then $|s-s_i|<2\delta$. 
Since $Q_t(s_i,a_{t,i})=v_t(s_i)$, the Lipschitz estimates above give 
\[ 
v_t(s)-Q_t(s,a_{t,i}) = v_t(s)-v_t(s_i) + Q_t(s_i,a_{t,i})-Q_t(s,a_{t,i})\le 4C_t\delta. 
\]
Averaging over the weights $\varphi_i(s)$ therefore yields 
\[ 
v_t(s) - \int_A Q_t(s,a)\pi_t^\delta(da\mid s) \le 4C_t\delta. 
\]
A backward induction over $t$ then gives 
\[ 
0\le \sup_{\widetilde\pi^k\in\Pi_\infty} \bar V^{k,s} \bigl( \widetilde\pi^k,\bar{\bm\pi};G^{\rm MF} \bigr) - \bar V^{k,s} \bigl( \pi^\delta,\bar{\bm\pi};G^{\rm MF} \bigr) \le 4\delta\sum_{t=0}^T C_t. 
\]
Finally, taking $\delta=C_S/L_\pi$ gives $\pi^\delta\in\Pi_{L_\pi}$ and hence
\[ 
0\le \sup_{\widetilde\pi^k\in\Pi_\infty} \bar V^{k,s} \bigl( \widetilde\pi^k,\bar{\bm\pi};G^{\rm MF} \bigr) - \sup_{\widetilde\pi^k\in\Pi_{L_\pi}} \bar V^{k,s} \bigl( \widetilde\pi^k,\bar{\bm\pi};G^{\rm MF} \bigr) \le \frac{C}{L_\pi}, 
\]
where $C:=4C_S\sum_{t=0}^T C_t$.
The corresponding bound on average exploitability follows directly by
integrating the preceding estimate with respect to $\mu_0^k$ and summing
over $k\in[K]$ with weights $\lambda_k$.

\subsection{Proof of Theorem \ref{thm:mf_approx}}
\label{proof.mf_approx}

$C$ denotes a finite positive constant which
may change from line to line and depends only on
$
T, 
W_{\max}, 
L_\pi,\ 
L_P, 
L_R,
E.
$
For notational simplicity, write
$$
\alpha_j:=\alpha_{N_j}(E),
\qquad
W_P:=\sum_{j=1}^K w_P^j,
\qquad
W_R:=\sum_{j=1}^K w_R^j.
$$
By definition of $W_{\max}$,
$
W_P\leq W_{\max},$ and $W_R\leq W_{\max}.$

We divide the proof into four steps.

\medskip

\noindent
{\bf Step 1. Independent non-identically distributed comparison system.}

Let
$$
\left(
L_t^1,\ldots,L_t^K
\right)_{t=0}^T
$$
be the state-action flow generated in
$G^{\mathrm{MF}}$ by the strategy profile
$\bar{\boldsymbol\pi}$.
For every $i\in I_k$, define a comparison process
$
(\bar X_t^i,\bar A_t^i)_{t=0}^T
$
as follows. Set
$
\bar X_0^i=s_0^i.
$
Given $\bar X_t^i$, sample
$$
\bar A_t^i
\sim
\bar\pi_t^k(\cdot\mid\bar X_t^i),
$$
and, for $t=0,\ldots,T-1$, sample
$$
\bar X_{t+1}^i
\sim
\bar P_t^k
\left(
\cdot
\mid
\bar X_t^i,\bar A_t^i,
L_t^1,\ldots,L_t^K
\right).
$$
The comparison processes are taken independent across
$i\in[N]$.
They are generally not identically distributed, since their
deterministic initial states $s_0^i$ may differ.
Let
$$
\nu_t^i
:=
\mathcal L(\bar X_t^i,\bar A_t^i).
$$
We claim that, for every $k$ and $t$,
\begin{equation}
\frac1{N_k}
\sum_{i\in I_k}\nu_t^i
=
L_t^k.
\label{eq:mf-average-law}
\end{equation}
Indeed, at time $0$,
$$
\frac1{N_k}
\sum_{i\in I_k}
\mathcal L(\bar X_0^i)
=
\frac1{N_k}
\sum_{i\in I_k}\delta_{s_0^i}
=
\mu_0^{k,N}.
$$
Applying the common policy $\bar\pi_0^k$ gives
$$
\frac1{N_k}
\sum_{i\in I_k}\nu_0^i
=
L_0^k.
$$
If the assertion holds at time $t$, then for every Borel set
$B\subset S$,
$$
\begin{aligned}
\frac1{N_k}
\sum_{i\in I_k}
\mathbb P(\bar X_{t+1}^i\in B)
&=
\frac1{N_k}
\sum_{i\in I_k}
\int_E
\bar P_t^k
\left(
B
\mid
s,a,L_t^1,\ldots,L_t^K
\right)
\nu_t^i(ds,da)
\\
&=
\int_E
\bar P_t^k
\left(
B
\mid
s,a,L_t^1,\ldots,L_t^K
\right)
L_t^k(ds,da) =
\mu_{t+1}^k(B).
\end{aligned}
$$
Applying $\bar\pi_{t+1}^k$ proves \eqref{eq:mf-average-law} at time $t+1$.
Define the empirical measure of the comparison agents by
$$
\bar L_t^{k,N}
:=
\frac1{N_k}
\sum_{i\in I_k}
\delta_{(\bar X_t^i,\bar A_t^i)}.
$$
By Lemma~\ref{lem:empirical-BL}, independence and \eqref{eq:mf-average-law} imply
\begin{equation}
\mathbb E
W_1
\left(
\bar L_t^{k,N},L_t^k
\right)
\leq
\alpha_k
\label{eq:mf-comparison-empirical}
\end{equation}
for every $t$ and $k$.
\noindent

{\bf Step 2. Propagation estimate for the lifted profile.}
Let
$
(X_t^i,A_t^i)_{t=0}^T
$
denote the state-action process in $G_N^{\mathrm{MF}}$ under the
lifted profile $\boldsymbol\pi$, and define
$$
L_t^{k,N}
:=
\frac1{N_k}
\sum_{i\in I_k}
\delta_{(X_t^i,A_t^i)}.
$$

We couple the interacting system and the comparison system introduced
in Step~1.

By Assumption~\ref{ass.regularkernel}, we may couple
$A_t^i$ and $\bar A_t^i$, conditionally on
$(X_t^i,\bar X_t^i)$, so that
\begin{equation}
    \mathbb E
\left[
|A_t^i-\bar A_t^i|
\mid
X_t^i,\bar X_t^i
\right]
\leq
L_\pi |X_t^i-\bar X_t^i|. \label{eq:mf-action-coupling}
\end{equation}
Set
$
c_\pi:=1+L_\pi.
$
It follows that
\begin{equation}
\mathbb E
d_E
\left(
(X_t^i,A_t^i),
(\bar X_t^i,\bar A_t^i)
\right)
\leq
c_\pi
\mathbb E|X_t^i-\bar X_t^i|.
\label{eq:mf-joint-coupling}
\end{equation}

For $k\in[K]$, define
$$
\Delta_t^k
:=
\frac1{N_k}
\sum_{i\in I_k}
\mathbb E|X_t^i-\bar X_t^i|,\qquad
D_t^k
:=
\mathbb E
W_1
\left(
L_t^{k,N},L_t^k
\right).
$$

Pairing the $i$-th atom of $L_t^{k,N}$ with the $i$-th atom of
$\bar L_t^{k,N}$ and using \eqref{eq:mf-comparison-empirical}--\eqref{eq:mf-joint-coupling} gives
\begin{equation}
D_t^k
\leq
c_\pi\Delta_t^k+\alpha_k.
\label{eq:mf-empirical-distance}    
\end{equation}
We next control the state discrepancy. By Assumptions~\ref{ass.transition_lipschitz}
and~\ref{ass.weightedLipschitz},
\begin{eqnarray}
&
W_1
\Big(
\bar P_t^k(
\cdot
\mid
X_t^i,A_t^i,
L_t^{1,N},\ldots,L_t^{K,N}),
\bar P_t^k(
\cdot
\mid
\bar X_t^i,\bar A_t^i,
L_t^1,\ldots,L_t^K)
\Big)
\\
&\leq
L_P
\left(
|X_t^i-\bar X_t^i|
+
|A_t^i-\bar A_t^i|
\right)
+\sum_{j=1}^K
w_P^j
W_1
\left(
L_t^{j,N},L_t^j
\right).
\label{eq:mf-transition-bound}
\end{eqnarray}
Here we have also used $W_1(\mu,\nu)
\leq
\|\mu-\nu\|_{\mathrm{TV}}$.

Using an optimal coupling of the two transition kernels, averaging
over $i\in I_k$, and using \eqref{eq:mf-action-coupling}, we obtain
\begin{equation}
\Delta_{t+1}^k
\leq
L_Pc_\pi\Delta_t^k
+
\sum_{j=1}^K w_P^jD_t^j.
\label{eq:mf-state-recursion}
\end{equation}
Combining \eqref{eq:mf-empirical-distance} and \eqref{eq:mf-state-recursion} yields
$$
\Delta_{t+1}^k
\leq
L_Pc_\pi\Delta_t^k
+
c_\pi
\sum_{j=1}^K
w_P^j\Delta_t^j
+\sum_{j=1}^K
w_P^j\alpha_j.
$$

Let
$
\Delta_t:=\max_{1\leq k\leq K}\Delta_t^k.
$ Since $\Delta_0=0$,
$$
\Delta_{t+1}
\leq
c_\pi(L_P+W_P)\Delta_t
+
\sum_{j=1}^K w_P^j\alpha_j.
$$
A finite-horizon discrete Gronwall argument therefore gives
\begin{equation}
\max_{0\leq t\leq T}\Delta_t
\leq
C\sum_{j=1}^K w_P^j\alpha_j.
\end{equation}
Substituting this into \eqref{eq:mf-empirical-distance} gives
\begin{equation}
D_t^k
\leq
\alpha_k
+
C
\sum_{j=1}^K w_P^j\alpha_j
\label{eq:mf-empirical-bound}
\end{equation}
uniformly in $t$ and $k$.
We also need a player-level estimate. For $i\in I_k$, let
$
\delta_t^i
:=
\mathbb E|X_t^i-\bar X_t^i|.
$
The same argument as above gives
$$
\delta_{t+1}^i
\leq
L_Pc_\pi\delta_t^i
+
\sum_{j=1}^K w_P^jD_t^j.
$$
Using \eqref{eq:mf-empirical-bound} and $\delta_0^i=0$ gives
\begin{equation}
\max_{0\leq t\leq T}
\max_{i\in[N]}
\delta_t^i
\leq
C
\sum_{j=1}^K w_P^j\alpha_j.
\label{eq:mf-player-state-bound}
\end{equation}
Let $\bar V^{k,s_0^i}
\left(
\bar\pi^k,\bar{\boldsymbol\pi};
G^{\mathrm{MF}}
\right)$
denote the value in the mean-field model of a population-$k$
representative player initialized at $s_0^i$ and following
$\bar\pi^k$. By Assumptions~\ref{ass.weightedLipschitz}
and~\ref{ass.reward_lipschitz},
$$
\begin{aligned}
&
\left|
V^i(\boldsymbol\pi;G_N^{\mathrm{MF}})
-
\bar V^{k,s_0^i}
\left(
\bar\pi^k,\bar{\boldsymbol\pi};
G^{\mathrm{MF}}
\right)
\right|
\leq
\sum_{t=0}^T
\left[
L_Rc_\pi\delta_t^i
+
\sum_{j=1}^K
w_R^jD_t^j
\right].
\end{aligned}
$$
Using \eqref{eq:mf-empirical-bound}--\eqref{eq:mf-player-state-bound},
\begin{equation}
\left|
V^i(\boldsymbol\pi;G_N^{\mathrm{MF}})
-
\bar V^{k,s_0^i}
\left(
\bar\pi^k,\bar{\boldsymbol\pi};
G^{\mathrm{MF}}
\right)
\right|
\leq
C\max\{\bar R_{\max},1\}
\sum_{j=1}^K
(w_R^j+w_P^j)\alpha_j.
\label{eq:mf-value-bound}
\end{equation}

\medskip

\noindent
{\bf Step 3. Unilateral deviations.}

Fix $i\in I_k$ and an arbitrary deviation
$\widetilde\pi^i\in\Pi$.
Let
$
(X_t^{\ell,(i)},A_t^{\ell,(i)})_{\ell\in[N]}
$
denote the process in $G_N^{\mathrm{MF}}$ under the profile
$
(\boldsymbol\pi^{-i},\widetilde\pi^i).
$
We now construct an independent comparison system
$$
(\bar X_t^{\ell,(i)},\bar A_t^{\ell,(i)})_{\ell\in[N]}
$$
in the deterministic mean-field environment
$(L_t^1,\ldots,L_t^K)$.
All agents $\ell\neq i$ use their original policies
$\bar\pi^{k(\ell)}$, whereas agent $i$ uses
$\widetilde\pi^i$. All comparison agents remain independent.
Let
$$
\bar L_t^{j,N,(i)}
:=
\frac1{N_j}
\sum_{\ell\in I_j}
\delta_{(\bar X_t^{\ell,(i)},\bar A_t^{\ell,(i)})}.
$$
For $j\neq k$, the average law of the agents in group $j$ is still
$L_t^j$. In group $k$, only one of the $N_k$ comparison agents has
been replaced. Hence, if
$$
\bar\nu_t^{k,(i)}
:=
\frac1{N_k}
\sum_{\ell\in I_k}
\mathcal L
(\bar X_t^{\ell,(i)},\bar A_t^{\ell,(i)}),
$$
then
$$
W_1
\left(
\bar\nu_t^{k,(i)},L_t^k
\right)
\leq
\frac{2}{N_k}.
$$
Applying Lemma~\ref{lem:empirical-BL} to the independent comparison agents and using the triangle inequality gives
\begin{equation}
\mathbb E
W_1
\left(
\bar L_t^{j,N,(i)},L_t^j
\right)
\leq
\alpha_j+
\frac{2{\bf 1}_{\{j=k\}}}{N_k}.
\end{equation}

Define
$$
D_t^{j,(i)}
:=
\mathbb E
W_1
\left(
L_t^{j,N,(i)},L_t^j
\right),\quad
{\rm  where}\quad
L_t^{j,N,(i)}
=
\frac1{N_j}
\sum_{\ell\in I_j}
\delta_{(X_t^{\ell,(i)},A_t^{\ell,(i)})}.
$$
Also define
$$
\Delta_t^{j,(i)}
:=
\frac1{N_j}
\sum_{\ell\in I_j}
\mathbb E
|X_t^{\ell,(i)}-\bar X_t^{\ell,(i)}|.
$$
Since every policy involved belongs to $\Pi$, Assumption~\ref{ass.regularkernel} applies
uniformly also to the deviating policy. Repeating the argument of
Step~2 gives
\begin{equation}
D_t^{j,(i)}
\leq
c_\pi\Delta_t^{j,(i)}
+
\alpha_j+
\frac{2{\bf 1}_{\{j=k\}}}{N_k},
\label{eq:mf-dev-empirical-distance}
\end{equation}
and
\begin{equation}
\Delta_{t+1}^{j,(i)}
\leq
L_Pc_\pi\Delta_t^{j,(i)}
+
\sum_{\ell=1}^K
w_P^\ell D_t^{\ell,(i)}.
\end{equation}
Therefore
\begin{equation}
\max_{0\leq t\leq T}
\max_{1\leq j\leq K}
\Delta_t^{j,(i)}
\leq
C
\left[
\sum_{\ell=1}^K
w_P^\ell\alpha_\ell
+
\frac{w_P^k}{N_k}
\right].
\label{eq:mf-dev-state-bound}
\end{equation}
Combining \eqref{eq:mf-dev-empirical-distance} and \eqref{eq:mf-dev-state-bound} gives
\begin{equation}
D_t^{j,(i)}
\leq
\alpha_j
+
\frac{2{\bf 1}_{\{j=k\}}}{N_k}
+
C
\left[
\sum_{\ell=1}^K
w_P^\ell\alpha_\ell
+
\frac{w_P^k}{N_k}
\right].
\label{eq:mf-dev-empirical-bound}
\end{equation}
For the deviating player itself, define
$
\delta_t^{i,(i)}
:=
\mathbb E
|X_t^{i,(i)}-\bar X_t^{i,(i)}|.
$
As in \eqref{eq:mf-transition-bound},
$$
\delta_{t+1}^{i,(i)}
\leq
L_Pc_\pi\delta_t^{i,(i)}
+
\sum_{j=1}^K
w_P^jD_t^{j,(i)}.
$$
Using \eqref{eq:mf-dev-empirical-bound},
\begin{equation}
\max_{0\leq t\leq T}
\delta_t^{i,(i)}
\leq
C
\left[
\sum_{j=1}^K w_P^j\alpha_j
+
\frac{w_P^k}{N_k}
\right].
\label{eq:mf-dev-player-state-bound} 
\end{equation}

Let
$
\bar V^{k,s_0^i}
\left(
\widetilde\pi^i,\bar{\boldsymbol\pi};
G^{\mathrm{MF}}
\right)
$
denote the value of a population-$k$ representative player,
initialized at $s_0^i$, using the deviation
$\widetilde\pi^i$ in the fixed mean-field environment generated by
$\bar{\boldsymbol\pi}$. By Assumptions~\ref{ass.weightedLipschitz}
and~\ref{ass.reward_lipschitz},
\eqref{eq:mf-dev-empirical-bound}--\eqref{eq:mf-dev-player-state-bound} imply
\begin{equation}
\begin{aligned}
&
\left|
V^i
\left(
\boldsymbol\pi^{-i},\widetilde\pi^i;
G_N^{\mathrm{MF}}
\right)
-
\bar V^{k,s_0^i}
\left(
\widetilde\pi^i,\bar{\boldsymbol\pi};
G^{\mathrm{MF}}
\right)
\right|
\\
&\leq
C\max\{\bar R_{\max},1\}
\left[
\sum_{j=1}^K
(w_R^j+w_P^j)\alpha_j
+
\frac1{N_k}
\right].
\end{aligned}
\end{equation}
The estimate is uniform over
$\widetilde\pi^i\in\Pi$.

Since
$$
\frac{1}{m}\le C\,\alpha_m(E),
\qquad m\ge1,
$$
the last term can be absorbed into $\alpha_k$. Thus
\begin{equation}
\begin{aligned}
&
\sup_{\widetilde\pi^i\in\Pi}
\left|
V^i
\left(
\boldsymbol\pi^{-i},\widetilde\pi^i;
G_N^{\mathrm{MF}}
\right)
-
\bar V^{k,s_0^i}
\left(
\widetilde\pi^i,\bar{\boldsymbol\pi};
G^{\mathrm{MF}}
\right)
\right|
\leq
C\max\{\bar R_{\max},1\}
\left[
\sum_{j=1}^K
(w_R^j+w_P^j)\alpha_j
+
\alpha_k
\right].
\end{aligned}
\label{eq:mf-dev-value-uniform}
\end{equation}

\medskip

\noindent
{\bf Step 4. Comparison of exploitabilities.}

For $i\in I_k$, define the state-conditioned mean-field
exploitability
$$
\operatorname{Expl}^{k,s_0^i}
\left(
\bar{\boldsymbol\pi};G^{\mathrm{MF}}
\right)
:=
\sup_{\widetilde\pi\in\Pi}
\left\{
\bar V^{k,s_0^i}
\left(
\widetilde\pi,\bar{\boldsymbol\pi};
G^{\mathrm{MF}}
\right)
-
\bar V^{k,s_0^i}
\left(
\bar\pi^k,\bar{\boldsymbol\pi};
G^{\mathrm{MF}}
\right)
\right\}.
$$

Using
$$
\sup_\eta f(\eta)-\sup_\eta g(\eta)
\leq
\sup_\eta|f(\eta)-g(\eta)|,
$$
together with \eqref{eq:mf-value-bound} and \eqref{eq:mf-dev-value-uniform}, we obtain
\begin{eqnarray}
\operatorname{Expl}_i
(\boldsymbol\pi;G_N^{\mathrm{MF}})
-
\operatorname{Expl}^{k,s_0^i}
(\bar{\boldsymbol\pi};G^{\mathrm{MF}})
\leq
C\max\{\bar R_{\max},1\}
\left[
\sum_{j=1}^K
(w_R^j+w_P^j)\alpha_j
+
\alpha_k
\right].
\label{eq:mf-exploitability-bound}
\end{eqnarray}
Average \eqref{eq:mf-exploitability-bound} over all players. Since
$$
\operatorname{AvgExpl}
(\boldsymbol\pi;G_N^{\mathrm{MF}})
=
\frac1N
\sum_{k=1}^K
\sum_{i\in I_k}
\operatorname{Expl}_i
(\boldsymbol\pi;G_N^{\mathrm{MF}}),
$$
and, by the definition of population-average exploitability in the
mean-field game,
$$
\operatorname{AvgExpl}
(\bar{\boldsymbol\pi};G^{\mathrm{MF}})
=
\frac1N
\sum_{k=1}^K
\sum_{i\in I_k}
\operatorname{Expl}^{k,s_0^i}
(\bar{\boldsymbol\pi};G^{\mathrm{MF}}),
$$
we obtain
$$
\begin{aligned}
&
\operatorname{AvgExpl}
(\boldsymbol\pi;G_N^{\mathrm{MF}})
-
\operatorname{AvgExpl}
(\bar{\boldsymbol\pi};G^{\mathrm{MF}})
\\
&\leq
C\max\{\bar R_{\max},1\}
\left[
\sum_{j=1}^K
(w_R^j+w_P^j)\alpha_{N_j}(E)
+
\sum_{k=1}^K
\frac{N_k}{N}\alpha_{N_k}(E)
\right].
\end{aligned}
$$
This is precisely the asserted bound.

\subsection{Proof of Theorem \ref{thm:perturbation}}\label{sec:proof-perturbation}
Throughout the proof, $C$ denotes a constant depending only on
$T,W_{\max},L_\pi,L_P,L_R$, whose value may change from line to line.

Fix an arbitrary strategy profile $\boldsymbol{\sigma}\in\Pi^N$, and couple the state-action processes $(X_t^i,A_t^i)$ and $(\widehat X_t^i,\widehat A_t^i)$ in $G$ and $G_N^{\mathrm{MF}}$ under the same strategy profile. Set
\[
\delta_t
:=
\max_{i\in[N]}
\mathbb E|X_t^i-\widehat X_t^i|.
\]
As in Step~2 of the proof of Theorem~\ref{thm:mf_approx},
Assumption~\ref{ass.regularkernel} and pairing corresponding atoms within
each population give
\[
\mathbb E
d_E\bigl(
(X_t^i,A_t^i),
(\widehat X_t^i,\widehat A_t^i)
\bigr)
\le
C\delta_t,
\]
and
\[
\mathbb E
W_1\bigl(
L_t^{k,N},
\widehat L_t^{k,N}
\bigr)
\le
C\delta_t.
\]

For $i\in I_k$, by the triangle inequality,
\[
\begin{aligned}
&W_1\!\left(
P_t^i(\cdot\mid \mathbf X_t,\mathbf A_t),
\widehat P_t^i(\cdot\mid
\widehat{\mathbf X}_t,\widehat{\mathbf A}_t)
\right)\\
&\quad\le
W_1\!\left(
P_t^i(\cdot\mid \mathbf X_t,\mathbf A_t),
\widehat P_t^i(\cdot\mid \mathbf X_t,\mathbf A_t)
\right)+
W_1\!\left(
\widehat P_t^i(\cdot\mid \mathbf X_t,\mathbf A_t),
\widehat P_t^i(\cdot\mid
\widehat{\mathbf X}_t,\widehat{\mathbf A}_t)
\right).
\end{aligned}
\]
The first term is bounded by $\epsilon_t^P$, while Assumptions 2 and 3,
together with the preceding coupling estimates and
$W_1\le \|\cdot\|_{\mathrm{TV}}$, bound the second term by $C\delta_t$.
Hence
\[
\delta_{t+1}\le C\delta_t+\epsilon_t^P.
\]
A finite-horizon discrete Gronwall estimate therefore yields
\[
\max_{0\le t\le T}\delta_t
\le
C\sum_{t=0}^{T-1}\epsilon_t^P.
\]

For $i\in[N]$, let
\[
r_t^i
:=
\sup_{\mathbf{s}\in S^N,\mathbf{a}\in A^N}
\left|
R_t^i(\mathbf{s},\mathbf{a})
-
\widehat R_t^i(\mathbf{s},\mathbf{a})
\right|.
\]
By Assumptions~\ref{ass.weightedLipschitz} and
\ref{ass.reward_lipschitz},
\[
\left|
V^i(\boldsymbol{\sigma};G)
-
V^i(\boldsymbol{\sigma};G_N^{\mathrm{MF}})
\right|
\le
C\left(
\sum_{t=0}^T r_t^i
+
\sum_{t=0}^{T-1}\epsilon_t^P
\right).
\]
The estimate is uniform over $\boldsymbol{\sigma}\in\Pi^N$.
As in Step~3 of the proof of Theorem~\ref{thm:mf_approx}, the same
estimate applies when player $i$ uses an arbitrary admissible deviation.
Hence
\[
\left|
\operatorname{Expl}^i(\boldsymbol{\pi};G)
-
\operatorname{Expl}^i(\boldsymbol{\pi};G_N^{\mathrm{MF}})
\right|
\le
C\left(
\sum_{t=0}^T r_t^i
+
\sum_{t=0}^{T-1}\epsilon_t^P
\right).
\]
Averaging over $i$ and using Cauchy--Schwarz,
\[
\frac1N\sum_{i=1}^N r_t^i
\le
\left(
\frac1N\sum_{i=1}^N(r_t^i)^2
\right)^{1/2}
=
\epsilon_t^R,
\]
gives the result.

\end{document}